\documentclass[11pt]{article}
\usepackage{amsfonts,a4wide}
\usepackage{amssymb}
\usepackage{graphicx}
\usepackage{mathtools}
\newtheorem{theorem}{Theorem}
 \newtheorem{corollary}[theorem]{Corollary}
 \newtheorem{lemma}[theorem]{Lemma}
 \newtheorem{proposition}[theorem]{Proposition}
 \newtheorem{definition}[theorem]{Definition}
\newtheorem{remark}[theorem]{Remark}
\newtheorem{example}[theorem]{Example}

\newcommand{\diag}{\mbox{\rm diag}}

\newcommand{\spa}{\mbox{\rm span}}

\DeclareMathOperator{\im}{im}

\newcommand{\Real}{\mathbb{R}}
\newcommand{\Comp}{\mathbb{C}}

\newcommand{\Field}{\mathbb{F}}

\newcommand{\eps}{\varepsilon}

\newcommand{\Q}[1]{\mathcal{Q}(#1)}

\newcommand{\set}[1]{\left\{#1\right\}}

\newcommand{\norm}[1]{\left\Vert#1\right\Vert}

\newcommand{\proof}{\par\noindent{\bf Proof}. \ignorespaces}
\newcommand{\eproof}{\space
    {\ \vbox{\hrule\hbox{\vrule height1.3ex\hskip0.8ex\vrule}\hrule}}\par}
 
\providecommand{\norm}[1]{\left\Vert#1\right\Vert}
\newcommand {\mycomment}[1]{} 
\newcommand {\mat}  [1] {\left[\begin{array}{#1}}
\newcommand {\rix}      {\end{array}\right]}

\newcommand {\ve}{\varepsilon}
\catcode`@=11     
\font\tenex=cmex10 
\newdimen\p@renwd
\setbox0=\hbox{\tenex B} \p@renwd=\wd0 
\def\bmat#1{\begingroup \m@th
  \setbox\z@\vbox{\def\cr{\crcr\noalign{\kern2\p@\global\let\cr\endline}}%
    \ialign{$##$\hfil\kern2\p@\kern\p@renwd&\thinspace\hfil$##$\hfil
      &&\quad\hfil$##$\hfil\crcr
      \omit\strut\hfil\crcr\noalign{\kern-\baselineskip}%
      #1\crcr\omit\strut\cr}}%
  \setbox\tw@\vbox{\unvcopy\z@\global\setbox\@ne\lastbox}%
  \setbox\tw@\hbox{\unhbox\@ne\unskip\global\setbox\@ne\lastbox}%
  \setbox\tw@\hbox{$\kern\wd\@ne\kern-\p@renwd\left[\kern-\wd\@ne
    \global\setbox\@ne\vbox{\box\@ne\kern2\p@}%
    \vcenter{\kern-\ht\@ne\unvbox\z@\kern-\baselineskip}\,\right]$}%
  \null\;\vbox{\kern\ht\@ne\box\tw@}\endgroup}
\catcode`@=12    
\newif\ifMDlatex
\catcode`@=11     
\def\MD@us#1{\csname#1style\endcsname}
\def\MD@uf#1{\csname#1font\endcsname}
\def\MD@t{text}
\def\MD@s{script}
\def\MD@ss{scriptscript}
\newdimen\MD@unit
\MD@unit\p@
\def\MD@changestyle#1{
  \relax\MD@unit0.1\fontdimen6\MD@uf{#1}0
  \everymath\expandafter{\the\everymath\MD@us{#1}}
}
\def\MD@dot{$\m@th\ldotp$}
\def\MD@palette#1{\mathchoice{#1\MD@t}{#1\MD@t}{#1\MD@s}{#1\MD@ss}}
\def\MD@ddots#1{{\MD@changestyle{#1}%
  \mkern1mu\raise7\MD@unit\vbox{\kern7\MD@unit\hbox{\MD@dot}}%
  \mkern2mu\raise4\MD@unit\hbox{\MD@dot}%
  \mkern2mu\raise \MD@unit\hbox{\MD@dot}\mkern1mu}}%
\def\MD@iddots#1{{\MD@changestyle{#1}%
  \mkern1mu\raise \MD@unit\hbox{\MD@dot}%
  \mkern2mu\raise4\MD@unit\hbox{\MD@dot}%
  \mkern2mu\raise7\MD@unit\vbox{\kern7\MD@unit\hbox{\MD@dot}}}}%
\def\MD@vdots#1{\vbox{\MD@changestyle{#1}%
    \baselineskip4\MD@unit\lineskiplimit\z@
    \kern6\MD@unit\hbox{\MD@dot}\hbox{\MD@dot}\hbox{\MD@dot}}}%
\ifMDlatex
  \DeclareRobustCommand\ddots{\mathinner{\MD@palette\MD@ddots}}%
  \DeclareRobustCommand\iddots{\mathinner{\MD@palette\MD@iddots}}%
  \DeclareRobustCommand\vdots{\mathinner{\MD@palette\MD@vdots}}%
\else
  \def\ddots{\mathinner{\MD@palette\MD@ddots}}%
  \def\iddots{\mathinner{\MD@palette\MD@iddots}}%
  \def\vdots{\mathinner{\MD@palette\MD@vdots}}%
\fi
\catcode`@=12    

\newcommand{\matp}[1]{\begin{bmatrix} #1 \end{bmatrix}}

\usepackage{color}

\begin{document}
\title
{Linear algebra properties of dissipative Hamiltonian descriptor systems}
\author{C. Mehl \footnotemark[3]~\footnotemark[1]
\and V. Mehrmann\footnotemark[3]~\footnotemark[1]
\and  M. Wojtylak \footnotemark[2]~\footnotemark[4]
}
\maketitle
\begin{abstract}

A wide class of matrix pencils connected with dissipative Hamiltonian
descriptor systems is investigated. In particular,  the following properties are shown: all eigenvalues are in the closed left half plane,
the nonzero finite eigenvalues on the imaginary axis are semisimple, the index is at most two, and
there are restrictions for the possible left and right minimal indices. For the case that the eigenvalue
zero is not semisimple, a structure-preserving method is presented that perturbs the given system into
a Lyapunov stable system.

\end{abstract}
\noindent

\noindent


\renewcommand{\thefootnote}{\fnsymbol{footnote}}
\footnotetext[3]{
Institut f\"ur Mathematik, Sekr.~MA 4-5, TU Berlin, Stra{\ss}e des 17.~Juni 136,
10623 Berlin, Germany.
\texttt{$\{$mehl,mehrmann$\}$@math.tu-berlin.de}.
}
\footnotetext[2]{Instytut Matematyki, Wydzia\l{} Matematyki i Informatyki,
Uniwersytet Jagiello\'nski, Krak\'ow, ul. \L ojasiewicza 6, 30-348 Krak\'ow, Poland
   \texttt{michal.wojtylak@uj.edu.pl}.}
\footnotetext[4]{   Supported by the Alexander von Humboldt Foundation with
  a  Scholarship for Experienced Scientists (carried out at TU Berlin) and by a Return Home Scholarschip.}
\footnotetext[1]{
Partially supported by the {\it Einstein Stiftung Berlin} through
the  Research Center {\sc Matheon} {\it Mathematics for key technologies}
in Berlin.
}
\renewcommand{\thefootnote}{\arabic{footnote}}

\section{Introduction}\label{sec:intro}

This paper deals with the linear algebra properties of matrix pencils that are associated with
linear time-invariant \emph{dissipative Hamiltonian descriptor systems} of the form
 \begin{equation}
E\dot x=\left(J-R\right)Qx. \label{diss_sys}
\end{equation}
These systems arise in energy based modeling of dynamical systems and are a special case of
so-called \emph{port-Hamiltonian
(pH) descriptor systems}, see, e.g., \cite{BeaMXZ17,Bre08,GolSBM03,JacZ12,OrtSMM01,Sch04,Sch13,Sch06}.
General pH descriptor systems have the form
\begin{eqnarray}
E\dot x&=&\left(J-R\right)Qx+(B-P)u,\nonumber \\
y&=& (B+P)^\star Qx + (S+N) u, \label{PHdef}
\end{eqnarray}
where for $(\Field,\star)\in\set{(\Comp,*),(\Real,\top)}$, we have $J,R\in\Field^{n,n}$, $E,Q\in\Field^{n,m}$,
$m\leq n$, and in addition $J^\star=-J$ and $R^\star=R\geq 0$, where we use the standard notation $W>0$ ($W\geq 0$)
for a real symmetric or complex Hermitian matrix $W$ to be positive (semi-)definite.

An important quantity of the system is the \emph{Hamiltonian} ${\mathcal H}(x)=\frac12 x^\star E^\star Qx$
which describes the distribution of internal energy among the energy storage elements of the system.
In order for this to make sense and to guarantee that $\mathcal H(x)$ is real, it is natural to assume that
\begin{equation}
E^\star Q=Q^\star E. \label{B1a}
\end{equation}
The matrix $J$ is the \emph{structure matrix} describing the energy flux among
energy storage elements within the system; $R$ is the \emph{dissipation matrix} describing energy
dissipation/loss in the system. The matrices $B\pm P\in\Field^{n,m}$ are \emph{port} matrices
describing the manner in which energy enters and exits the system, and $S+N$, with $S=S^\star\in \Field^{m,m}$
and $N=-N^\star \in \Field^{m,m}$, describe the direct \emph{feed-through} from input to output.
In most practical examples the system has the extra property that $Q^\star E$ is positive semidefinite,
in which case we speak of a pHDAE \emph{with nonnegative Hamiltonian}.
%
%
%
%
%
%

To illustrate the broad variety of application areas where pH descriptor systems play a role,
we review the following examples from the literature.
\begin{example}\label{ex:circuit}{ \rm
A simple RLC network, see, e.g., \cite{BeaMXZ17,Dai89,Fre11}, can be modeled by a DAE of the form
\begin{eqnarray}\label{eq:RLC_network_1}
\underbrace{\mat{ccc} G_c C G_c^\top & 0 & 0\\ 0 & L & 0 \\ 0 & 0 & 0 \rix}_{:=E}
\mat{c}\dot{V} \\ \dot{I_l}\\ \dot{I_v} \rix=
\underbrace{\mat{ccc} -G_r R_r^{-1}G_r^\top & -G_l & -G_v \\ G_l^\top & 0 & 0\\ G_v^\top & 0 & 0 \rix}_{:=J-R}
\mat{c}V \\ I_l\\ I_v \rix,
\end{eqnarray}
where $L > 0$, $C > 0$, $R_r> 0$ are real symmetric, $G_v$ is of full column rank, and
the other matrices have appropriate dimensions. (For a detailed physical interpretation of the matrices
and variables we refer to \cite{BeaMXZ17}.) Here, $J$ and $-R$ are defined to be the skew-symmetric and symmetric parts,
respectively, of the matrix on the right hand side of~\eqref{eq:RLC_network_1}. We see that
in this example we have a dissipative Hamiltonian descriptor system of the form \eqref{diss_sys} with the matrix $Q$ being the identity,
$E=E^T\geq 0$, $J=-J^T$, and $R\geq 0$.
}
\end{example}

\begin{example} \label{ex:stokes} {\rm Space discretization of the Stokes or Oseen equation in fluid dynamics,
see, e.g., \cite{EmmM13}, leads to a dissipative Hamiltonian system with pencil
\[
\lambda E-(J-R)=\lambda \matp{M&0\\0 &0} - \matp{A & B\\ -B^\top& 0},\quad J=\mat{cc}0&B\\ -B^\top&0\rix,\quad
R=\mat{cc}A&0\\ 0&0\rix,
\]
with $A$ being a positive semidefinite discretization of the negative Laplace operator, $B$ being a discretized gradient,
and $M$ being a positive definite mass matrix, so $R\geq 0$, $Q=I$, and  $E=E^\top\geq 0$. Note that the matrix $B$ may
be rank deficient, in which case the pencil is singular.}
\end{example}

\begin{example} \label{ex:gas}{\rm Space discretization of the Euler equation describing the flow in a gas network \cite{EggKLMM17} leads to a dissipative Hamiltonian descriptor system with
\[
E=\matp{M_1& 0 &0\\ 0 & M_2 & 0\\ 0 & 0 & 0},\quad J=\matp{0& -G &0\\ G^\top  & 0 & K^\top\\ 0 & -K & 0},\quad R=\matp{0& 0 &0\\ 0 & D & 0\\ 0 & 0 & 0},\quad Q=I,
\]
where $M_1,M_2,D$ are positive definite and the matrices $N$ and $\matp{G^\top & K^\top}$ have full column rank.
Thus, again we have $R\geq 0$ and a nonnegative Hamiltonian, since $Q^\top E=E^\top Q\geq 0$.
}
\end{example}

\begin{example}\rm
Classical second order representations of linear damped mechanical systems, see, e.g., \cite{TisM01,Ves11},  have the form
\begin{equation}\label{eq:19.10.17a}
M\ddot x +D\dot x + Kx= f,
\end{equation}
where $M,D,K\in\mathbb R^{n,n}$ are symmetric positive semidefinite matrices and represent the mass, damping,
and stiffness matrices, respectively, while $f$ is an exiting force that is assumed to be sufficiently
often differentiable. First order formulation leads to  a system associated with the matrix pencil
\begin{equation}\label{eq:19.10.17}
\lambda\underbrace{\mat{cc}M&0\\ 0&I\rix}_{=:E}
-\left(\vphantom{\mat{cc} 0&I\\ -I&0\rix}\right.
\underbrace{\mat{cc} 0&I\\ -I&0\rix}_{=:J}
- \underbrace{\mat{cc}D&0\\ 0&0\rix}_{=:R}
\left.\vphantom{\mat{cc}D&0\\ 0&0\rix}\right)
     \underbrace{\mat{cc}I &0\\ 0&K\rix}_{=:Q},
\end{equation}
where we have $E^\top Q=Q^\top E\geq 0$ and $R\geq 0$. As a particular case, freely vibrating strongly damped
systems of the form $M\ddot x +sD\dot x + Kx=0$ were investigated in \cite{Tas15} for $s\to\infty$.
In the limit case, the first order formulation of the equivalent system $\frac{1}{s}M\ddot x+D\dot x+\frac{1}{s}Kx=0$
leads to a pencil as in~\eqref{eq:19.10.17} with $M=K=0$. Observe that this pencil is singular if
$D$ is singular.
\end{example}

\begin{example}\label{ex:mech}{\rm  Consider a mechanical system as in~\eqref{eq:19.10.17a}
with a position constraint $Gx=0$. Replacing this constraint
by its derivative $G\dot x=0$  and incorporating it in the dynamical equation using a Lagrange multiplier $y$,
the system  can be written  as
\[
 \matp{M & 0 \\ 0 & 0}\matp{\ddot x \\ \ddot  y}+\matp{D & 0\\ 0 & 0}\matp{\dot x\\ \dot y}+\matp{K & G^\top \\ G & 0 }\matp{x\\ y}=\matp{f\\0},
\]
which in first order formulation gives a system associated with the pencil
\[
\lambda\underbrace{\matp{M &&&\\ &0&&\\ &&I& \\ &&& I }}_{=:E}-
\left(\vphantom{\matp{ &&I&\\ &&&I\\ -I&&& \\ &I&&  }}\right.
\underbrace{\matp{ &&I&\\ &&&I\\ -I&&& \\ &I&&  }}_{=:J}
- \underbrace{\matp{D\\ &0&&\\ &&0& \\ &&& 0 }}_{=:R}
\left.\vphantom{\matp{D\\ &0&&\\ &&0& \\ &&& 0 }}\right)
     \underbrace{\matp{I &&&\\ &I&&\\ &&K&G^\top \\ &&G& 0 }}_{=:Q}.
\]
Here, both $E$ and $Q$ are symmetric, commuting square matrices and $R$ is positive semidefinite,
but in contrast to the previous examples $E^\top Q$ is not positive semidefinite in this case.
Note also that often the system is modeled with redundant constraints which implies that $G$ is rank deficient and then the pencil is singular.
}
\end{example}

Similar examples arise in many other applications, see e.g. \cite{GraMQSW16}, where the singular pencil
case typically arises as a limiting case or when redundant system modeling is used.

To understand the linear algebra properties of pH descriptor systems, we analyze matrix pencils
of the form
\begin{equation}\label{phpencil}
P(\lambda)=\lambda E-LQ\in\Field^{n,m}[\lambda],\text{ with }
 L=J-R,\quad J=-J^\star,\quad R=R^\star\geq 0,
\end{equation}
where $\Field^{n,m}[\lambda]$ denotes the set of $n\times m$ matrix polynomials over  $\Field$.
So far, mainly the cases of $Q=I$ and $Q>0$ were considered in the literature, but here
we will deal with the general case, allowing even that the pencil is rectangular and thus singular.
If $Q=I$ and $R=0$, then $P(\lambda)$ is an \emph{even pencil}, i.e., the coefficient associated with $\lambda$
is Hermitian/symmetric and the second coefficient is skew-Hermitian/skew-symmetric.
Canonical forms for $P(\lambda)$ in these
cases are well known, see, e.g., \cite{Tho76}. If, in addition, $E$ is positive semidefinite,
then it easily follows from these forms that all finite eigenvalues are on the imaginary axis and
semisimple. Then replacing the condition $R=0$ by $R\geq 0$ has the effect that now
all finite eigenvalues are in the closed left-half complex plane, and those on the imaginary axis are still semisimple, see, e.g., \cite{MehMS16}.

At first sight, it does not seems obvious that the structure of the pencil $P(\lambda)$ in~\eqref{phpencil}
preserves any of these restrictions in the spectrum. However, it turns out that the assumption of both
$R$ and $Q^\star E$ being positive semidefinite leads to surprisingly rich linear algebraic properties
of $P(\lambda)$.

The remainder of the paper is organized as follows.
Section \ref{s:prel} reviews the Kronecker canonical form and introduces the concept of
regular  deflating subspaces of singular matrix pencils. It appears that the properties of the pencil $\lambda E-Q$ are crucial for understanding the pencil $P(\lambda)$. Hence, we we first analyze the pencil $\lambda E-Q$ in Section \ref{sec:propEQ}. In the regular case a canonical form
is derived  in Propositions~\ref{prop:EQreg} and \ref{EQrectangle}. In the singular case we present a
condensed form (Theorem \ref{thm:condensed})  and show what are the possible left and right minimal indices
(Corollary \ref{cor:indices}). Sections~\ref{sec:ph} and~\ref{sec:phproof} discuss the stability properties
and the Kronecker canonical form of the pencil  $P(\lambda)=\lambda E - LQ$.
In particular, we show in Theorem \ref{thm:singind} that if the pencil $\lambda E-Q$ does not have
left minimal indices larger than zero,
then the pencil $P(\lambda)$ has all finite eigenvalues in the closed left half-plane,
and  all eigenvalues on the imaginary axis except zero are semisimple. Furthermore, the index of the pencil (see Section~\ref{s:prel}) is at most two and its right minimal indices (if there are any) are not larger than one. If, in addition,
the pencil $\lambda E-Q$ is regular, then the left minimal indices of $P(\lambda)$ are all zero
(if there are any). Furthermore, we present many examples that illustrate which properties are lost if some of
the assumptions on the pencil $\lambda E-Q$ are weakened.
Finally, we discuss in Section~\ref{sec5} structure-preserving perturbations of $P(\lambda)$ that make
the eigenvalue zero semisimple, but keep the Jordan structure of all other eigenvalues invariant.


\section{Preliminaries}\label{s:prel}

For general matrix pencils $\lambda E-A$ the structural properties are characterized via the Kronecker canonical form.
Recall the following result for complex matrix pencils \cite{Gan59a}.
\begin{theorem}\label{th:kcf}
Let $E,A\in {\mathbb C}^{n,m}$. Then there exist nonsingular matrices
$S\in {\mathbb C}^{n,n}$ and $T\in {\mathbb C}^{m,m}$ such that
\begin{equation}\label{kcf}
S(\lambda E-A)T=\diag({\cal L}_{\epsilon_1},\ldots,{\cal L}_{\epsilon_p},
{\cal L}^\top_{\eta_1},\ldots,{\cal L}^\top_{\eta_q},
{\cal J}_{\rho_1}^{\lambda_1},\ldots,{\cal J}_{\rho_r}^{\lambda_r},{\cal N}_{\sigma_1},\ldots,
{\cal N}_{\sigma_s}),
\end{equation}
where the block entries have the following properties:
\begin{enumerate}
\item[\rm (i)]
Every entry ${\cal L}_{\epsilon_j}$ is a bidiagonal block of size
${\epsilon_j}\times ({\epsilon_j+1})$, $\epsilon_j\in{\mathbb N}_0$,
of the form
\[
\lambda\left[\begin{array}{cccc}
1&0\\&\ddots&\ddots\\&&1&0
\end{array}\right]-\left[\begin{array}{cccc}
0&1\\&\ddots&\ddots\\&&0&1
\end{array}\right].
\]
\item[\rm (ii)]
Every entry ${\cal L}^\top_{\eta_j}$ is a bidiagonal block of size
$({\eta_j+1})\times {\eta_j}$, $\eta_j\in{\mathbb N}_0$,
of the form
\[
\lambda\left[\begin{array}{ccc}
1\\0&\ddots\\&\ddots&1\\&&0
\end{array}\right]-
\left[\begin{array}{ccc}
0\\1&\ddots\\&\ddots&0\\&&1
\end{array}\right].
\]
\item[\rm (iii)]
Every entry ${\cal J}_{\rho_j}^{\lambda_j}$ is a Jordan block of size
${\rho_j}\times{\rho_j}$, $\rho_j\in{\mathbb N}$, $\lambda_j\in{\mathbb C}$,
of the form
\[
\lambda\left[\begin{array}{cccc}
1\\&\ddots\\&&\ddots\\&&&1
\end{array}\right]-
\left[\begin{array}{cccc}
\lambda_j&1\\&\ddots&\ddots\\&&\ddots&1\\&&&\lambda_j
\end{array}\right].
\]
\item[\rm (iv)]
Every entry ${\cal N}_{\sigma_j}$ is a nilpotent block of size
${\sigma_j}\times {\sigma_j}$, $\sigma_j\in{\mathbb N}$,
of the form
\[
\lambda\left[\begin{array}{cccc}
0&1\\&\ddots&\ddots\\&&\ddots&1\\&&&0
\end{array}\right]-
\left[\begin{array}{cccc}
1\\&\ddots\\&&\ddots\\&&&1
\end{array}\right].
\]
\end{enumerate}
The Kronecker canonical form is unique up to permutation of the blocks.
\end{theorem}

For real matrices there exists a real Kronecker canonical form which is obtained under real transformation
matrices $S,T$. Here, the  blocks ${\cal J}_{\rho_j}^{\lambda_j}$ with $\lambda_j\in\Comp\setminus\Real$ are in real Jordan canonical form instead, but the other
blocks have the same structure as in the complex case.
%

The sizes $\eta_j$, and $\epsilon_i$ of the rectangular blocks
are called the \emph{left and right minimal indices} of $\lambda E-A$, respectively. Furthermore,
a value $\lambda_0\in\mathbb C$ is called a (finite) eigenvalue of $\lambda E-A$ if
\[
\operatorname{rank}(\lambda_0E-A)<\operatorname{nrank}(\lambda E_A):=\min_{\alpha\in\mathbb C}
\operatorname{rank}(\alpha E-A),
\]
where $\operatorname{nrank}(\lambda E_A)$ is called the \emph{normal rank} of $\lambda E-A$. Furthermore,
$\lambda_0=\infty$ is said to be an eigenvalue of $\lambda E-A$ if zero is an eigenvalue of $\lambda A-E$.
It is obvious from Theorem~\ref{th:kcf} that the blocks $\mathcal J_{\rho_j}$ as in (iii) correspond to
finite eigenvalues of $\lambda E-A$, whereas blocks $\mathcal N_{\sigma_j}$ as in (iv) correspond to the
eigenvalue $\infty$. The sum of all sizes of blocks that are associated with a fixed eigenvalue
$\lambda_0\in\mathbb C\cup\{\infty\}$ is called the \emph{algebraic multiplicity} of $\lambda_0$.
The size of the largest block ${\cal N}_{\sigma_j}$ is
called the \emph{index} $\nu$ of the pencil $\lambda E-A$, where, by convention,  $\nu=0$ if $E$ is invertible.

The matrix pencil $\lambda E-A\in\mathbb F^{n,m}[\lambda]$ is called \emph{regular} if $n=m$ and
$\operatorname{det}(\lambda_0 E-A)\neq 0$ for some $\lambda_0 \in \mathbb C$,
otherwise it is called \emph{singular}. A pencil is singular if and only if it has blocks of at least one of the types ${\cal L}_{\eps_j}$ or
${\cal L}^\top_{\eta_j}$ in the Kronecker canonical form.

While eigenvalues of singular pencils are well-defined, the same is not true for eigenvectors.
For example, the pencil
\begin{equation}\label{exdefl}
\lambda E-A=\mat{ccc}\lambda-1&0&0\\ 0&\lambda&-1\rix={\cal J}_1^1\oplus {\cal L}_1
\end{equation}
has the eigenvalue $1$ with algebraic multiplicity one, but any vector $x$ of the form
$x=\mat{ccc}\alpha&\beta&\beta\rix^\top$ with $\alpha,\beta\in\mathbb F$ satisfies $\lambda_0 Ex=Ax$.
Hence, the equation $\lambda_0 Ex=Ax$ itself is not suitable (at least for the purpose of the current paper)
as a generalization of the notion \emph{eigenvector} in the singular case. We therefore
introduce the following concept.
\begin{definition}\label{def:regdefsub}\rm
Let $E,A\in\mathbb F^{n,m}$ and let $\lambda_0\in\mathbb C$ be an eigenvalue of $P(\lambda)=\lambda E-A$ with
algebraic multiplicity $k$.
\begin{enumerate}
\item[(i)] A $k$-dimensional subspace $\mathcal X$ of $\mathbb C^m$ is called a \emph{(right) regular deflating
subspace of $P(\lambda)$ associated with $\lambda_0$} if there exists nonsingular matrices $X\in\Comp^{m,m}$ and $Y\in\mathbb C^{n,n}$
such that the first $k$ columns of $X$ span $\mathcal X$ and such that
\begin{equation}\label{eq:regdefsub}
Y(\lambda E-A)X=\mat{cc}R(\lambda)&0\\ 0&\widetilde P(\lambda)\rix,
\end{equation}
where $R(\lambda)$ is a regular $k\times k$ pencil that has only the eigenvalue $\lambda_0$, and where
$\widetilde P(\lambda)$ does not have $\lambda_0$ as an eigenvalue.
\item[(ii)]
A vector $x\in\mathbb C^m\setminus\{0\}$ is called a \emph{regular eigenvector} of $P(\lambda)$ associated
with $\lambda_0$, if $\lambda_0Ex=Ax$ and if there exists a regular deflating subspace associated with
$\lambda_0$ that contains $x$.
\end{enumerate}
Regular  deflating subspaces and regular eigenvectors associated with the eigenvalue $\infty$ are defined in
the obvious way by considering the reversed pencil $\lambda A-E$ and the eigenvalue zero.
\end{definition}
Clearly, for each eigenvalue of a singular pencil there exists a regular eigenvector and a regular deflating subspace.
For the pencil in \eqref{exdefl} each subspace of the form ${\cal X}=\spa(\mat{ccc}\alpha&\beta&\beta\rix^\top)$
with some $\alpha,\beta\in\mathbb F$, $\alpha\neq 0$ is a regular deflating subspace. Indeed, with
\[
X=\mat{ccc} \alpha & 0 &0 \\ \beta & 1 & 0\\ \beta &0 & 1\rix  ,\quad Y=\mat{cc} 1/\alpha & 0 \\ -\beta/\alpha & 1 \rix
\]
one has $Y(\lambda E-A) X=\lambda E-A$.
In particular, this shows that a regular deflating subspace associated with an eigenvalue is not unique if the pencil is singular, which is in contrast with the regular case. Furthermore, it is straightforward to
show that neither the subspace
\[
\spa(\mat{ccc}0 & 1& 1\rix^\top)\quad\mbox{nor}\quad
\spa(\mat{ccc}0 & 1& 1\rix^\top, \mat{ccc}1 & 0& 0\rix^\top)
\]
is a regular deflating subspace. We have the following key property of regular deflating subspaces.
%
\begin{lemma}\label{lem:regdefsub}
Let $E,A\in\mathbb F^{n,m}$, let $\lambda_0\in\mathbb C\cup\{\infty\}$ be an eigenvalue of
$\lambda E-A$, and let $\mathcal X$ be an associated regular deflating subspace. Then
\[
{\cal X}\cap \spa\left( \bigcup_{\lambda\in(\Comp\cup\set\infty)\setminus\{\lambda_0\}}\!\!\!\!\ker (\lambda E- A)\right)=\set0,
\]
with the usual convention $\infty E-A:=E$.
\end{lemma}
\proof
Indeed, let $X,Y$ be as in Definition~\ref{def:regdefsub} such
that~\eqref{eq:regdefsub} is satisfied. Furthermore assume that $y\in\ker (\lambda_1 E- A)$ for some
$\lambda_1\in(\Comp\cup\set\infty)\setminus\{\lambda_0\}$ and let
$X^{-1}y=\mat{cc}y_{1}^\top  & y_{2}^\top\rix^\top$ be split conformably with \eqref{eq:regdefsub}. First we show that $y_1=0$. Indeed,
\[
0= Y (\lambda_1 E-A) X( X^{-1}y)=   \mat{c}R(\lambda_1) y_1 \\ *\rix.
\]
In particular, $y\in\ker R(\lambda_1)$. But the pencil $R(\lambda)$ is regular and has only one eigenvalue $\lambda_0\neq\lambda_1$. Therefore $y_1=0$. In other words, $PX^{-1}y=0$, where $P$ is the canonical projection onto the first $k$ coordinates.
Assume now that $x=\sum_{i=1}^l z_i$, with $x\in\cal X$ and $z_i\in\ker(\lambda_i E-A)$, $\lambda_i\neq \lambda_0$ for $i=1,\dots,l$ and some $l>0$.
Then $PX^{-1}x= \sum_{i=1}^l Pz_i=0$. But since $x$ is a vector in $\mathcal X$, there exists a vector $v\in\Field^k$
such that $x=X\mat{cc}v^\top&0\rix^\top$, i.e., $PX^{-1}x=X^{-1}x$. Hence $x=0$, which finishes the proof.
%
%
\eproof
Let us mention two other concepts associated  with eigenvectors and deflating subspaces in the case of singular pencils that are present in the literature.
In \cite{Van83},  \emph{reducing subspaces} were introduced
to replace the concept of deflating subspaces in the singular case. This concept has the advantage that
reducing subspaces associated with an eigenvalue are now uniquely defined. However, they are explicitly
allowed to contain vectors \emph{from the singular part} and thus a result similar to  Lemma~\ref{lem:regdefsub} cannot be obtained.
For the second concept, that only refers to eigenvectors, we first discuss a different characterization of regular eigenvectors as follows.
Let $\mathcal N_r(P)$ denote the right nullspace of $P(\lambda)$ interpreted as a matrix over the field of rational functions (thus $\mathcal N_r(P)$ is a subspace of the space $\mathbb C^m(\lambda)$ of vectors of
length $m$ with rational functions as entries) and let
\[
\mathcal N_r(P)\big|_{\lambda_0}:=\big\{z\in\mathbb C^m\,\big|\,z=v(\lambda_0)\mbox{ for some }
v(\lambda)\in\mathcal N_r(P)\big\}
\]
denote the \emph{right nullspace of $P(\lambda)$ evaluated at $\lambda_0$}, see \cite{Fro16}.
Then it can be shown that $x\in\mathbb C^m$ is a regular eigenvector of $P(\lambda)$ associated with $\lambda_0$  if and only
if $\lambda_0Ex=Ax$ and $x\not\in\mathcal N_r(P)\big|_{\lambda_0}$.
This construction allows to define \emph{eigenspaces} uniquely, e.g., as being
orthogonal to $\mathcal N_r(P)\big|_{\lambda_0}$, cf. \cite{Fro16}.
However, this orthogonality property is not preserved under equivalence transformations for  pencils.

\section{Properties of the pencil $\lambda E-Q$}\label{sec:propEQ}
As it was announced in the introduction, in this section we discuss the properties of a (possibly rectangular) pencil
$\lambda E-Q\in\Field^{m,n}[\lambda]$ satisfying the condition~\eqref{B1a}.
We will also frequently assume that, in addition, one (or both) of the conditions
 \begin{equation}
EQ^\star=QE^\star,\label{B1b}
\end{equation}
and
\begin{equation}
E^\star Q\geq 0\label{B1c}
\end{equation}
%
%
are satisfied. Although \eqref{B1b} does not have a direct physical interpretation,
it is satisfied by most classical examples (see the introduction) and together with
\eqref{B1a} it forms a very strong assumption on the structure of the pencil $\lambda E-Q$. Condition \eqref{B1c}
has a minor impact on the properties of $\lambda E-Q$, but strongly influences the properties
of $P(\lambda)$, as will be seen in subsequent sections. We will not discuss pencils $\lambda E-Q$ for which only
\eqref{B1b} is  assumed, because results for these pencils can easily be obtained from the corresponding
results on pencils satisfying~\eqref{B1a}  by considering the pencil $\lambda E^\star-Q^\star$.

Clearly, when $E$ and $Q$ are square Hermitian (or real symmetric matrices) that commute (as is the case
in some of our applications), then both $E^\star Q=Q^\star E$ and $EQ^\star=QE^\star$ are satisfied.
It is clear that for any unitary $U\in\Field^{m,m}$ and invertible $X\in\Field^{n,n}$ a transformation of
the form $(E,Q)\mapsto (UEX,UQX)$ preserves the conditions~\eqref{B1a} and~\eqref{B1c}. The following
result yields a canonical form for regular pencils under this transformation.
%
\begin{proposition}\label{prop:EQreg}
Let $\lambda E-Q\in\mathbb F^{n,n}[\lambda]$ be a regular pencil. Then we have $E^\star Q=Q^\star E$ if and
only if there exist a unitary (in the case $\mathbb F=\mathbb R$ real orthogonal) matrix $U\in\mathbb F^{n,n}$ and an
invertible $X\in\mathbb F^{n,n}$ such that $UEX$ and $UQX$ are diagonal and real, and satisfy
\[
(UEX)^2+(UQX)^2=I_{n}.
\]
Moreover, $U$ and $X$ can be chosen such that one of the matrices $UEX$ or $UQX$ is nonnegative.
If, in addition, we have $E^\star Q\geq 0$, then  $U$ and $X$ can be chosen such that $UEX$ and $UQX$ are both nonnegative.
\end{proposition}
\proof
Suppose that $E^\star Q=Q^\star E$ holds. Then the regularity of the pencil $\lambda E-Q$ implies that the
columns of  the matrix
\[
L:=\matp{E\\Q}
\]
are linearly independent and they span an $n$-dimensional \emph{Lagrangian subspace}, i.e., a subspace
$\mathcal{U}\subset \Field^{2n}$ such that $$u^T \matp{0 & I_n\\ -I_n & 0}v=0$$
for all $u,v\in\mathcal{U}$. Orthonormalizing the columns of $L$, i.e., carrying out a factorization,
$L=ZT$ with $T\in \Field^{n,n}$ invertible, and
\[
Z=\matp{Z_1 \\ Z_2},\ Z_1,Z_2 \in \Field^{n,n}
\]
satisfying $Z^\star Z=I_n$, we obtain that the columns of $Z=LT^{-1}$ span the same Lagrangian subspace
as the columns of $L$ and according to the CS decomposition (see, e.g., \cite[Theorem 2.1]{PaiV81} for the
complex case and \cite[Lemma 2.1]{MehP16} for the real case), there exists $U,V\in\mathbb F^{n,n}$ unitary and such that
\[
\matp{U &0\\ 0&U} \matp{E\\Q}TV =\matp{C\\S},
\]
with $C$ and $S$ being diagonal such that $C^2+S^2=I_n$. The claim follows now by taking $X:=TV$.
The converse implication is straightforward.

Clearly, by premultiplying $U$ with a diagonal signature matrix that switches the signs of diagonal
entries if necessary then yields the desired nonnegativity of one of the matrices $UEX$ or $UQX$.
If, in addition, we have $E^\star Q\geq 0$, then it follows that
\[
(UEX) (UQX)=(UEX)^\star (UQX)=X^\star E^\star Q X\geq 0
\]
and
thus nonnegativity of one of the diagonal  matrices $UEX$ or $UQX$ implies  nonnegativity of the other.
\eproof
A direct corollary of Proposition~\ref{prop:EQreg} gives a characterization of the index of $\lambda E-Q$.
\begin{corollary}\label{cor:index}
Let $E,Q\in\mathbb F^{n,n}$ be such that $E^\star Q=Q^\star E$ and that the pencil $\lambda E-Q$ is regular.
Then it is of index at most one.
\end{corollary}
\proof
The proof follows directly from Proposition~\ref{prop:EQreg}, since if both $E$ and $Q$ are transformed
to be diagonal, then the regularity implies that the nilpotent part in the Kronecker canonical form
is diagonal and hence the index is zero if $E$ is nonsingular and it is one otherwise.
\eproof
The next result shows that the assumption of regularity in Proposition~\ref{prop:EQreg} can be dropped,
if both conditions~\eqref{B1a} and~\eqref{B1b} are satisfied. Note that the rectangular case $m\neq n$ is included. By a diagonal matrix $D=\mat{c}d_{ij}\rix$ we then mean a matrix with all entries $d_{ij}$
equal to zero unless $i=j$.
\begin{proposition}\label{EQrectangle}
Let $E,Q\in\Field^{m,n}$ be such that $E^\star Q=Q^\star E$ and $EQ^\star=QE^\star$. Then there exist
unitary (real orthogonal), matrices $U\in\mathbb F^{m,m},V\in\mathbb F^{n,n}$ such that $UEV$ and $UQV$ are
 diagonal and real.
Furthermore, $U$ and $V$ can be chosen such that one of the matrices $UEV$ and $UQV$ is nonnegative, and
if in addition $E^\star Q\geq 0$ holds, then $U,V$ can be chosen such that $UEV$ and $UQV$ are both
nonnegative.
\end{proposition}
\proof
The proof follows from the fact that if both conditions $E^\star Q=Q^\star E$ and $EQ^\star=QE^\star$ hold then this is equivalent to the fact that the four matrices $E^\star Q$, $Q^\star E$, $Q^\star Q$ and $E^\star E$ are all
normal and commuting, as well as the four matrices $EQ^\star$, $QE^\star$, $QQ^\star$ and $EE^\star$. Then the
statement follows immediately from \cite[Corollary 5]{maehara2011simultaneous} in the complex case
and \cite[Corollary 9]{maehara2011simultaneous} in the real case.
\eproof
Surprisingly, the statement of Proposition~\ref{prop:EQreg} is in general not true for singular pencils $\lambda E-Q$
satisfying only $E^\star Q=Q^\star E$ as the following example shows.
\begin{example}\label{ex:sp} {\rm Let
\[
E=\matp{1&0\\1&1\\0&0}=S\matp{1&0\\0&1\\0&0}, \quad Q=\matp{0&1\\1&0\\1&1}=S\matp{0&0\\1&0\\0&1},\
\mbox{with}\ S=\matp{1&0 & 1\\1&1&0\\0&0&1}.
\]
Then $E^\star Q$ is symmetric as
\[
E^\star Q=\matp{1 & 1 \\ 1&0},
\]
but by construction $\lambda E-Q$ has one singular block with left minimal index $2$
in its Kronecker canonical form and thus cannot be equivalent to a pencil in diagonal form.
}
\end{example}

Although a direct generalization of Proposition~\ref{prop:EQreg} to the singular case is not possible,
it turns out that there is a restriction in the possible Kronecker structure for pencils
satisfying~\eqref{B1a} and~\eqref{B1c}.
\begin{proposition}\label{prop:minind}
Let $E,Q\in\mathbb F^{m,n}$ be such that $E^\star Q=Q^\star E$. Then all right minimal indices of $\lambda E-Q$ are zero
(if there are any).
\end{proposition}
\proof
Let $\lambda E-Q= S(\lambda E_0 -Q_0)T$, where $\lambda E_0-Q_0$ is in Kronecker canonical form.
Clearly, with $E^\star Q$ also $E_0^\star S^\star S Q_0$ is Hermitian.
Suppose that $\lambda E-Q$ has a right minimal index which is nonzero. Then $\lambda E_0-Q_0$ contains a block
of the form
\begin{equation}\label{block}
\lambda\matp{1 &0 & \\ &\ddots&\ddots&\\ &&1&0}+\matp{0 &1 & \\ &\ddots&\ddots&\\ &&0&1}
\end{equation}
of size $k\times (k+1)$ with $k>0$. Without loss of generality let this block be the first block in the
Kronecker canonical form $\lambda E_0-Q_0$. Then the $k$th and $(k+1)$st standard basis vectors $x:=e_k$
and $y:=e_{k+1}$ satisfy $Q_0y=E_0x\neq 0$ and $E_0y=0$. This, however, implies that
\[
y^\star Q_0 ^\star S^\star S E_0 x=0<y^\star E_0^\star S^\star SQ_0x
\]
contradicting the fact that $E_0^\star S^\star S Q_0$ is Hermitian.
\eproof
Clearly, by applying Proposition~\ref{prop:minind} to the pencil $\lambda E^\star- Q^\star$ we get that  $EQ^\star=QE^\star$ implies  all left minimal indices of $\lambda E-Q$ being zero.
Another immediate consequence of Proposition~\ref{prop:minind} is the following corollary.
\begin{corollary}\label{cor:eleven}
Let $E,Q\in\mathbb F^{n,n}$ be such that $E^\star Q=Q^\star E$. If $\lambda E-Q$ is singular, then
$E$ and $Q$ have a common nullspace.
\end{corollary}
\proof
This follows immediately from Proposition~\ref{prop:minind}, the fact that a square singular
pencil must have both left and right minimal indices, and the fact that a right minimal index zero
corresponds to a common nullspace.
\eproof
With the help  of these intermediate results, we are now able to derive a Kronecker-like condensed form for pencils $\lambda E-Q$ satisfying $E^\star Q=Q^\star E$.
\begin{theorem}\label{thm:condensed}
Let $E,Q\in\mathbb F^{m,n}$ be such that $E^\star Q=Q^\star E$. Then there exist a unitary (real orthogonal) matrix
$U\in\mathbb F^{m,m}$ and a nonsingular matrix $X\in\mathbb F^{n,n}$ such that
\begin{equation}\label{cform}
UEX=\mat{ccc}E_{11} & E_{12} & 0\\ 0 & E_{22}&0\rix\quad\mbox{and}\quad
UQX=\mat{ccc}Q_{11} & Q_{12} & 0\\ 0 & Q_{22}&0\rix
\end{equation}
and the following statements hold.
\begin{itemize}
\item[\rm (i)] $E_{11},Q_{11}\in\mathbb F^{n_1,n_1}$ are diagonal, real and satisfy $E_{11}^2+Q_{11}^2=I_{n_1}$.
Furthermore, $U$ and $X$ can be chosen such that one of the matrices $E_{11}$ and $Q_{11}$ is nonnegative and if, in addition, $E^\star Q\geq 0$, then $U,X$ can be chosen such that both $E_{11}$ and $Q_{11}$ are
nonnegative.
\item[\rm (ii)] $Q_{22},E_{22}\in\mathbb F^{m_2,n_2}$, with $m_2=n-n_1> n_2$ $($or $m_2=n_2=0)$, are such that $\lambda E_{22}-Q_{22}$
is singular having left singular blocks (blocks of the form $\mathcal L^\top_{\eta_j}$) only;
\item[\rm (iii)] $\operatorname{Ker}(E_{22})\subseteq \operatorname{Ker}(E_{12})$ and
$\operatorname{Ker}(Q_{22})\subseteq \operatorname{Ker}(Q_{12})$.
\item[\rm (iv)]The  pencil $\lambda E-Q$ has the same Kronecker canonical form as
\begin{equation}\label{2pencils}
\lambda\mat{ccc}E_{11} & 0 & 0\\ 0 & E_{22}& 0\rix-\mat{ccc}Q_{11} & 0 & 0 \\ 0 & Q_{22} & 0\rix.
\end{equation}
\end{itemize}
In particular,  the regular part of $\lambda E-Q$ is given by $\lambda E_{11}-Q_{11}$,
 the left minimal indices of $\lambda E-Q$ coincide with those of $\lambda E_{22}-Q_{22}$ and
 $\lambda E-Q$ has exactly $n-n_1-n_2$ right minimal indices which are all zero.
\end{theorem}
\proof
 Let $S,T$ be nonsingular matrices such that
\begin{equation}\label{kcf_pencil}
S(\lambda E-Q)T=\mat{cc}\mathcal R(\lambda)&0\\ 0&\mathcal S(\lambda)\rix
\end{equation}
is in Kronecker canonical form, where $\mathcal R(\lambda)\in\mathbb F^{n_1,n_1}$ and
$\mathcal S(\lambda)\in\mathbb F^{m_1,n-n_1}$ are the regular part and singular part, respectively.
Without loss of generality, let the singular blocks in $\mathcal S(\lambda)$ be ordered in such a way
that the left singular blocks (of the form $\mathcal L^\top_{\eta_j}$) come first. Then, since by
Proposition~\ref{prop:minind} all right minimal indices of $\lambda E-Q$ are zero, we obtain
that $\mathcal S(\lambda)$ has the form
\[
\mathcal S(\lambda)=\mat{cc}\lambda E_0-Q_0&0\rix,
\]
where $\lambda E_0-Q_0\in\mathbb F^{m_2,n_2}$ contains all singular blocks ${\cal L}^\top_{\eta_j}$.
This implies that $m_2>n_2$ unless there are no blocks ${\cal L}_{\eta_j}^\top$, i.e., $m_2=n_2=0$.
The zero block is of size $n-n_1-n_2$, i.e.,  the pencil $\lambda E-Q$ has $n-n_1-n_2$
right minimal indices all of which are zero. Let
\[
S^{-1}=V\mat{cc}R_{11}&R_{12}\\ 0&R_{22}\rix,\quad \text{with} \quad V^\star V=I,
 \]
be a $QR$ decomposition of $S^{-1}$ partitioned conformably with (\ref{kcf_pencil}). Then, with
$\lambda \widetilde E_{11}-\widetilde Q_{11}:=R_{11}\mathcal R(\lambda)$ one has
\[
V^\star (\lambda E-Q)T= \mat{ccc}\lambda \widetilde E_{11}-\widetilde Q_{11} &R_{12}(\lambda E_0-Q_0) &0\\
0&R_{22}(\lambda E_0-Q_0) &0\rix.
\]
Furthermore,  $E^\star Q=Q^\star E$ implies that
$\widetilde E_{11}^\star\widetilde Q_{11}=\widetilde Q_{11}^\star\widetilde E_{11}$. Since
$R$, and thus also $R_{11}$ and $R_{22}$ are invertible, it follows that $\lambda\widetilde E_{11}-\widetilde Q_{11}$
is regular. Thus, by Proposition~\ref{prop:EQreg} there exists a unitary matrix $U_1\in\mathbb F^{n_1,n_1}$ and
a nonsingular matrix $X_1\in\mathbb F^{n_1,n_1}$ such that $E_{11}:=U_1 \widetilde E_{11}X_1$ and
$Q_{11}:=U_1\widetilde Q_{11}X_1$ are real diagonal satisfying $E_{11}^2+Q_{11}^2=I_{n_1}$, with at least one of the matrices $E_1,Q_1$ being nonnegative.
Setting  $U:=\diag(U_1,I_{m_2})V^\star $, $X:=T\diag(X_1,I_{n-n_1})$, we   obtain
\begin{equation}\label{E0Q022}
U (\lambda E-Q)X= \mat{ccc}\lambda U_1(\widetilde  E_{11}-\widetilde Q_{11})X_1 & U_1 R_{12}(\lambda E_0-Q_0) &0\\
0&R_{22}(\lambda E_0-Q_0) &0\rix,
\end{equation}
which is the desired form \eqref{cform}.
We also have already shown the first statement of (i) and statement (ii).
To finish the proof of (i), note that if  $E^\star Q$ is positive semidefinite, then
$\widetilde E_{11}^\star\widetilde Q_{11}$ is positive definite as well, as a principle submatrix of
a matrix congruent to $E^\star Q$. Hence, one may apply Proposition~\ref{prop:EQreg}.

 Invertibility of $R_{22}$ together with \eqref{E0Q022} implies that statement (iii) holds.
To see (iv) note that the equivalence transformation
$$
\mat{cc} I & 0 \\  -U_1 R_{12}R_{22}^{-1}  & I \rix U(\lambda E- Q) X
$$
transforms $\lambda E-Q$ into the pencil in \eqref{2pencils}, see also again equation \eqref{E0Q022}.
\eproof
So far we have not seen any restrictions on the left minimal indices of pencils
$\lambda E-Q$ satisfying $E^\star Q=Q^\star E$. In fact, the left minimal indices may be arbitrary, even under the
additional assumption $E^\star Q\geq 0$ as long as the
corresponding singular blocks are of suitable dimensions to be part of the block $\lambda E_{22}-Q_{22}$ in the condensed form of
Theorem~\ref{thm:condensed}. To show this, it would be helpful to obtain a canonical form for pencils $\lambda E-Q$
satisfying $Q^\star E=E^\star Q\geq 0$. The problem in obtaining such a form lies in the fact that
left-multiplication with a non-unitary invertible matrix $S$ does not preserve the property of $Q^\star E$
being Hermitian, let alone being Hermitian positive semi-definite. In particular, if $\lambda E_0-Q_0=\mathcal L_{\eta}^\top$
is a Kronecker block associated with a left minimal index $\eta>0$, then it is straightforward
to check that $Q_0^\star E_0$ is not a Hermitian matrix. We will therefore investigate the structure of matrices that make the pencil $\lambda E_0-Q_0$ in Kronecker canonical form Hermitian by left-multiplication.
\begin{proposition}\label{canlemma}
Let $S\in\Field^{l,k}$, and let $E_1=SE_0$, $Q_1=SQ_0$, where
\begin{equation}\label{block2}
\lambda E_0-Q_0=\mathcal L^\top_{k-1}=\lambda\matp{1 & &\\ 0 &\ddots& \\ &\ddots&1\\ && 0 }-\matp{0 && \\ 1 &\ddots&\\ &\ddots&0\\ && 1}\in\Comp^{k,k-1}[\lambda]
\end{equation}
Then the following statements hold.
\begin{itemize}
\item[\rm (i)] We have $E_1^\star Q_1=Q_1^\star E_1$ if and only if $H=S^\star S$ is a positive definite real Hankel matrix,
i.e.,
\begin{equation}\label{eq:hankel}
H=\mat{cccc}  a_0 & a_1 & \dots &a_{k-1} \\ a_1 &\iddots&\iddots&a_{k} \\  \vdots &\iddots&\iddots& \vdots\\
a_{k-1} & a_k & \dots & a_{2k-1}  \rix\in\Real^{k,k},\quad H>0.
\end{equation}
\item[\rm (ii)] We have $E_1^\star Q_1=Q_1^\star E_1\geq 0$ (or, equivalently, $E_1^\star Q_1=Q_1^\star E_1>0$)
if and only if $H=S^\star S$ is a  positive definite real Hankel matrix such that its submatrix $\tilde H$
obtained by deleting the first row and last column is also positive definite, i.e.,
\begin{equation}\label{eq:hankel2}
\tilde H=\mat{cccc}  a_1 & a_2 & \dots &a_{k-1} \\ a_2 &\iddots&\iddots&a_{k} \\  \vdots &\iddots&\iddots& \vdots\\ a_{k-1} & a_{k+1} & \dots & a_{2k-2}  \rix > 0.
\end{equation}
\end{itemize}
\end{proposition}
\proof
(i) Suppose that $H=S^\star S$ is real and has the form given in~\eqref{eq:hankel}. Then
%
\begin{equation}\label{tildeh}
E_1^\star Q_1= E_0^\star H Q_0=\tilde H= Q_0^\star H E_0=Q_1^\star E_1,
\end{equation}
where $\tilde H$ has the form as in~\eqref{eq:hankel2} (but need not be positive definite). Since this matrix
is a real Hankel matrix and thus Hermitian, it follows that $E_1^\star Q_1$ is Hermitian.

Conversely, assume that $E_1^\star Q_1=Q_1^\star E_1$. Then it follows by comparison of elements in the
equation $E_0^\star H Q_0=Q_0^\star H E_0$ that $H$ is a Hermitian Hankel matrix, or equivalently,
a real Hankel matrix.

(ii) follows immediately from (\ref{tildeh}). Note that for $E_1^\star Q_1$ positive semi-definiteness
is equivalent to positive definiteness, since $E_1^\star Q_1$ has full rank.
\eproof
%
%
%
%

Using these preparations, we now obtain the following result showing that a pencil $\lambda E-Q$
satisfying $Q^\star E=E^\star Q\geq0$ may have arbitrary left minimal indices as long as they fit
the dimension restrictions.
\begin{corollary}\label{cor:indices}
Let $n,m\in\mathbb N$,  and let $q,\eta_1,\dots, \eta_{q}$ be arbitrary non-negative integers satisfying
\[
n-q\leq m, \quad \eta_1+\cdots +\eta_{q}\leq n-q.
\]
Then there exists a pencil $\lambda E-Q\in\Field^{m,n}[\lambda]$
satisfying $Q^\star E=E^\star Q\geq 0$ with precisely $q$ left minimal indices $\eta_1,\dots, \eta_{q}$.
\end{corollary}
%
\proof
First observe that for each $k\geq 1$ there exists a positive definite Hankel matrix $H$ with $\tilde H>0$,
where $\tilde H$ is defined as in~\eqref{eq:hankel2}, when $H$ is as in~\eqref{eq:hankel}.
Indeed, let $\xi_1>\cdots >\xi_k>0$ and let
\[
V=\mat {ccc}  \xi_1 & \dots & \xi_1^{n} \\  \vdots &  & \vdots\\   \xi_{n} & \dots & \xi_{n}^{n}  \rix   \in\Real^{k,k}
\]
be a Vandermonde matrix. Observe that
\begin{equation}\label{eq:Hvander}
 H=V^\top V=\mat{c} \sum\limits_{l=1}^n \xi_l^{i+j}\rix_{i,j=1}^k
\end{equation}
is a positive definite Hankel matrix. Furthermore,  $\tilde H>0$  due to
\[
x^\star \tilde H x= \sum_{i,j=1}^{k-1} x_i\bar x_j \sum_{l=1}^k  \xi_l^{i+j+1} =\sum_{l=1}^k  \xi_l
\left |\sum_{i=1}^{k-1} x_i \xi_l^i\right|^2  >0,\quad\mbox{for } x\in\Field^{k-1}\setminus\{0\}.
\]
%
%
%
Taking $S=H^{1/2}$ (or $S$ to be the Cholesky factor of $H$) we get by Proposition \ref{canlemma} (ii) a pencil $\lambda E_1-Q_1$ satisfying $E_1^\star Q_1\geq0$ with left minimal index $k-1$.
To obtain a pencil with the desired minimal indices
it is enough to take the direct sum of pencils constructed as above for the minimal indices $\eta_1,\dots,\eta_q$
and if necessary adding to the direct sum an appropriate pencil having no left minimal indices.
(It is possible to find such a pencil, because its size is
 $(n-\eta-q)\times (m-\eta)$ with $\eta=\eta_1+\cdots +\eta_q$ which means that the number of rows does not exceed the number of columns by the assumption $n-q\leq m$.)
\eproof
The natural equivalence relation for the set of pencils $\lambda E-Q$ satisfying $E^\star Q=Q^\star E$ is
$(\lambda E_1-Q_1)\sim (\lambda E_2-Q_2)$ if and only if there exist unitary $U$ and invertible $X$ such that
$U(\lambda E_1-Q_1)X= \lambda E_2-Q_2$.
So far we have established the canonical form and  the complete set of invariants  in the regular case
(Proposition~\ref{prop:EQreg}) and in the singular case if $EQ^\star =QE^\star$ is additionally assumed
(Proposition~\ref{EQrectangle}). We have also established all the possible Kronecker canonical forms
in the general case (Theorem~\ref{thm:condensed} and Corollary~\ref{cor:indices}). However, in establishing
the complete set of invariants of the equivalence relation $\sim$ in the general case, we have met the following difficulty.
\begin{remark}{\rm
For any $k>1$ there is a continuum of nonequivalent (in the sense of the relation $\sim$ defined above)
linear pencils
$\lambda E-Q\in\Field^{k,k-1}[\lambda]$ with $E^\star Q=Q^\star E\geq0$ having precisely one left minimal
index $k-1$. This can be seen in the following way.
Suppose that two pencils $\lambda S_1E_0 -Q_0$ and $\lambda S_2 E_0-Q_0$ as in Proposition~\ref{canlemma}
(ii) are equivalent, i.e., we have
$U^\star (\lambda S_1E_0 -Q_0)X=\lambda S_2 E_0-Q_0$ with unitary $U$ and invertible $X$. Then by
Proposition~\ref{canlemma}~(ii) the matrix $H_i:=S_i^\star S_i\in\Real^{k,k}$ ($i=1,2$) is a Hankel matrix.
The relation $\sim$ implies that
\[
X^\star E_0^\star H_1 E_0 X=  E_0^\star H_2 E_0 ,\quad X^\star Q_0^\star H_1 Q_0 X=  Q_0^\star H_2 Q_0
\]
and thus the Hermitian pencils $\lambda  E_0^\star H_1 E_0-  Q_0^\star H_1 Q_0$ and
$\lambda  E_0^\star H_2 E_0-  Q_0^\star H_2 Q_0$ are congruent.

Now let   $\xi_1>\dots>\xi_n>0$  and let $H$ be the positive definite Hankel matrix defined as
in \eqref{eq:Hvander}.
Observe that there  clearly is a continuum of mutually non-congruent pencils of the form
\[
\lambda  E_0^\star H E_0-  Q_0^\star H Q_0 =
\lambda E_0^\star\mat{c} \sum_{l=1}^n \xi_l^{i+j}\rix_{i,j=1}^{k-1}E_0 -
Q_0^\star\mat{c} \sum_{l=1}^n \xi_l^{i+j}\rix_{i,j=2}^kQ_0.
\]
One can see this easily, e.g., by investigating the $\xi_1$ dependence of eigenvalues, keeping
$\xi_2,\dots,\xi_n$ fixed. But this immediately implies that under the equivalence relation $\sim$ there is a continuum of nonequivalent pencils $\lambda E-Q$ with canonical form $\lambda E_0-Q_0$ that satisfy $E^\star Q=Q^\star E\geq 0$.
}
\end{remark}

In this section, we have considered the possible canonical forms for pencils $\lambda E-Q$ satisfying one
or more of the conditions~\eqref{B1a}--\eqref{B1c}. In the next section, we will discuss the influence
of these conditions on the pencil $P(\lambda)$ from~\eqref{phpencil}.

\section{Properties of the pencil $P(\lambda)=\lambda E-(J-R)Q$}\label{sec:ph}
In this section, we investigate the properties of the pencil $P(\lambda)=\lambda E-(J-R)Q$ from~\eqref{phpencil}.
For brevity, we will sometimes use the notation $L:=J-R$ and thus write our pencil in the form
\begin{equation}\label{phpencilwithl}
P(\lambda)=\lambda E-LQ,\quad E,Q\in\mathbb F^{n,m},\quad L\in\mathbb F^{m,m},\quad L+L^\star\leq 0.
\end{equation}
Note that the decomposition $L=J-R$, with $J=-J^\star$, $R=R^\star$ is unique. Hence, the condition $L+L^\star\leq 0$ is equivalent to the condition $R\geq 0$.
\begin{proposition}\label{prop:regular}
Suppose that the pencil $P(\lambda)$ in~\eqref{phpencilwithl} is regular with $E^\star Q=Q^\star E$. Then
the pencil $\lambda E-Q$ is regular as well.
\end{proposition}
\proof
Suppose that $\lambda E-Q$ is singular. Then by Corollary~\ref{cor:eleven} the matrices $E$ and $Q$ have a
common right nullspace which implies that $E$ and $LQ$ have a common right nullspace as well, contradicting the
hypothesis that $P(\lambda)$ is regular.
\eproof
In the regular case and with $Q>0$, it is well-known, see e.g., \cite{GilMS18,MehMS16}, that all eigenvalues
of the pencil $P(\lambda)$ in~\eqref{phpencilwithl} are in the closed left half complex plane and those
on the imaginary axis are semisimple. We will now investigate if these properties remain valid under
the weaker condition $E^\star Q=Q^\star E\geq 0$. Clearly, without this condition, we may have eigenvalues in
the right half complex plane as is shown by the scalar example
\[
E=\mat{c}-1\rix, \quad J=\mat{c}0\rix\quad R=\mat{c}1\rix,
\quad Q=\mat{c}1\rix,\quad L=J-R=\mat{c}1\rix,
\]
where the corresponding pencil $\lambda E-LQ=1-\lambda$
has the eigenvalue $\lambda_0=1$. But even if the condition $E^\star Q=Q^\star E\geq 0$ is satisfied we
may have eigenvalues in the right half plane if the
pencil $\lambda E-Q$ happens to be singular.
\begin{example}\label{ex:rhp}\rm
Consider the matrices
\[
E=\mat{cc}1&0\\ 0&0\rix,\quad Q=\mat{cc}0&0\\ a&0\rix,\quad L=\mat{cc}0&1\\ -1&0\rix,
\]
where $a>0$. Then we have $E^*Q=Q^*E=0\geq 0$ and $L+L^\star=0\leq 0$, but the pencil
\[
P(\lambda)=\lambda\mat{cc}1&0\\ 0&0\rix-\mat{cc}a&0\\ 0&0\rix
\]
has the eigenvalue $a>0$ in the right half plane.
\end{example}

We highlight that in the previous example, the pencil $\lambda E-Q$ had a left minimal index larger than zero.
(We recall that all right minimal indices are automatically zero by Theorem~\ref{thm:condensed}, but
the left minimal indices may be arbitrary up to dimensional restrictions by Corollary~\ref{cor:indices}.)
In fact, it turns out that under the additional hypothesis that all left minimal indices of $\lambda E-Q$
are zero, the pencil $P(\lambda)$ has all eigenvalues in the closed left half complex plane, and those
on the imaginary axis are semisimple with the possible exception of non semi-simplicity of the eigenvalues zero
and infinity.
With these preparations we are now able to formulate and prove our main result that does not only
give information on the finite eigenvalues, but also restricts the possible index and minimal indices
of a pencil as in~\eqref{phpencilwithl}.
\begin{theorem}\label{thm:singind}
Let $E,Q\in \mathbb F^{n,m}$ satisfy $E^\star Q=Q^\star E\geq 0$ and let all left minimal indices of
$\lambda E-Q$ be equal to zero (if there are any). Furthermore, let
$L\in\mathbb F^{m,m}$ be such that $R:=-\frac{1}{2}(L+L^\star)\geq 0$. Then the following
statements hold for the pencil $P(\lambda)=\lambda E-LQ$.
\begin{enumerate}
\item[\rm (i)] If $\lambda_0\in\mathbb C$ is an eigenvalue of $P(\lambda)$ then $\operatorname{Re}(\lambda_0)\leq 0$.
\item[\rm (ii)] If $\omega\in\mathbb R\setminus\{0\}$ and $\lambda_0=i\omega$ is an eigenvalue of $P(\lambda)$, then
$\lambda_0$ is semisimple. Moreover, if the columns of $V\in\mathbb C^{m,k}$ form a basis of a regular deflating
subspace of $P(\lambda)$ associated with $\lambda_0$, then $RQV=0$.
\item[\rm (iii)] The index of $P(\lambda)$ is at most two.
\item[\rm (iv)] All right minimal indices of $P(\lambda)$ are at most one (if there are any).
\item[\rm (v)] If in addition $\lambda E-Q$ is regular, then all left minimal indices of $P(\lambda)$ are zero
(if there are any).
\end{enumerate}
\end{theorem}

Before we give the proof of Theorem~\ref{thm:singind} in the next section, we discuss a few
consequences and examples.
\begin{remark}{\rm
Part (i) and (ii) of Theorem~\ref{thm:singind} generalize \cite[Lemma 3.1]{MehMS16} in which these results
were proved for the case $E=I$ and $Q>0$, and they also generalize \cite[Theorem 2]{GilMS18} where the regular
case was covered under the additional assumption that for any eigenvector $x$ associated with a nonzero
eigenvalue the condition $Q^\star Ex\neq 0$ is satisfied, but no necessary or sufficient conditions were
given when this is the case.
Theorem~\ref{thm:singind} now states that this extra assumption is not needed in the regular case
(in fact, the proof of Theorem~\ref{thm:singind} shows that this condition is automatically satisfied),
and it extends the assertion to the singular case under the additional assumption
that the left minimal indices of the pencil $\lambda E-Q$ are all zero.
}
\end{remark}

\begin{example}\label{rem:nonsimple0}{\rm
We highlight that Theorem~\ref{thm:singind} gives no information about the semi-simplicity of
the eigenvalue zero of the pencil $P(\lambda)$. Indeed, consider the matrices
\[
C=\mat{cc}1&0\\ 0&1\rix,\;J=\mat{cc}0&-1\\ 1&0\rix,\;D=\mat{cc}1&0\\ 0&0\rix,\;R=\mat{cc}0&0\\ 0&0\rix.
\]
Then with $E=C$ and $Q=D$ the pencil
\[
P(\lambda)=\lambda E-(J-R)Q=\lambda\mat{cc}1&0\\ 0&1\rix-\mat{cc}0&0\\ 1&0\rix
\]
satisfies $Q^\star E=D\geq 0$ and
has the eigenvalue zero with algebraic multiplicity two, but geometric multiplicity one.
(We will see in Section~\ref{sec5} that sizes of Jordan blocks associated with the eigenvalue zero are also not restricted to be at most two, even under the additional hypothesis that $\lambda E-Q$ is regular.)
}
\end{example}

\begin{example}{\rm
With the matrices from Example~\ref{rem:nonsimple0}, and $E=D$ and $Q=C$ the pencil
\[
P(\lambda)=\lambda E-(J-R)Q=\lambda\mat{cc}1&0\\ 0&0\rix-\mat{cc}0&-1\\ 1&0\rix
\]
satisfies again $Q^\star E=D\geq 0$ and
has the eigenvalue $\infty$ with algebraic multiplicity two, but geometric multiplicity one.
This shows that the case that the index of the pencil is two may indeed occur even when
$\lambda E-Q$ is regular.
}
\end{example}

\begin{example}{\rm
The case that the pencil $P(\lambda)$ as in Theorem~\ref{thm:singind} has a right minimal index
equal to one may indeed occur even in the case when $\lambda E-Q$ is regular. For an example, consider the matrices
\[
E=\mat{cc}1&0\\ 0&0\rix,\quad Q=\mat{cc}0&0\\ 0&1\rix,\quad L=\mat{cc}-1&1\\ -1&0\rix,\quad
LQ=\mat{cc}0&1\\ 0&0\rix.
\]
Then $\lambda E-Q$ is regular, $E^*Q=Q^*E=0$, and $L+L^*\leq 0$. The pencil
$P(\lambda)=\lambda E-LQ$ has one right minimal index equal to one and one left minimal index equal to zero.
}
\end{example}

\begin{example}\label{rem:ind}{\rm
If the pencil $\lambda E-Q$ in Theorem~\ref{thm:singind} is not regular, then not much can be said about
the left minimal indices of the pencil $P(\lambda)$. In fact, they can be arbitrarily large as the following
example shows. Let
\[
E=Q=\mat{cc}I_{n-1}&0\\ 0&0\rix,\quad J=\mat{ccccc}0&-1&0&\dots&0\\1&0&-1&\ddots&\vdots\\
0&\ddots&\ddots&\ddots&0\\ \vdots&\ddots&\ddots&\ddots&-1\\ 0&\dots&0&1&0\rix\in\mathbb R^{n,n},\quad
R=0\in\mathbb R^{n,n}
\]
and $L=J-R=J$. Then $Q^\star E=E^\star Q=Q^2\geq 0$ and the pencil $\lambda E-Q$ has one left and right
minimal index both being equal to zero. Furthermore, we obtain
\[
P(\lambda)=\lambda E-LQ=\mat{ccccc}\lambda&1&0&\dots&0\\-1&\lambda&\ddots&\ddots&0\\
0&\ddots&\ddots&1&0\\ \vdots&\ddots&\ddots&\lambda&0\\ 0&\dots&0&-1&0\rix.
\]
Obviously this pencil has one right minimal index which is zero. Since its normal rank is $n-1$,
it has exactly one left minimal index. Observing that $\operatorname{rank}E=\operatorname{rank}P(\lambda_0)=n-1$
for all $\lambda_0\in\mathbb C$, it follows that $P(\lambda)$ does not have finite or infinite eigenvalues,
and hence the left minimal index must be $n-1$ which is the largest possible size for a left minimal
index of an $n\times n$ singular matrix pencil.
}
\end{example}

As an application of Theorem~\ref{thm:singind}, we obtain the following statement on the Kronecker structure
of quadratic matrix polynomials that correspond to damped mechanical systems.
\begin{corollary}
Let $S(\lambda):=\lambda^2M+\lambda D+K\in\mathbb F^{n,n}[\lambda]$ be a quadratic matrix polynomial such that
$M,D,K$ are all Hermitian and positive semidefinite matrices. Then the following statements
hold.
\begin{enumerate}
\item[\rm (i)] All eigenvalues of $S(\lambda)$ are in the closed left half complex plane and all finite nonzero
eigenvalues on the imaginary axis are semisimple.
\item[\rm (ii)] The possible length of Jordan chains of $S(\lambda)$ associated with either with the eigenvalue $\infty$ or with the eigenvalue zero is at most two.
\item[\rm (iii)] All left and all right minimal indices of $S(\lambda)$ are zero (if there are any).
\end{enumerate}
\end{corollary}
\proof
The proof of (i) and the statement in (ii) on the length of the Jordan chains associated with the eigenvalue
$\infty$ follows immediately from the fact that the companion linearization
\[
P(\lambda)=\lambda\mat{cc}M&0\\ 0&I\rix-\mat{cc}D&K\\ -I&0\rix=\lambda\mat{cc}M&0\\ 0&I\rix-\left(\mat{cc}0&I\\ -I&0\rix
-\mat{cc}D&0\\ 0&0\rix\right)\mat{cc}I&0\\ 0&K\rix
\]
satisfies the hypothesis of Theorem~\ref{thm:singind} with
\[
E=\mat{cc}M&0\\ 0&I\rix,\quad J=\mat{cc}0&I\\ -I&0\rix,\quad R=\mat{cc}D&0\\ 0&0\rix,\quad Q=\mat{cc}I&0\\ 0&K\rix.
\]
The remaining statement of (ii) then follows by applying the already proved part of (ii)
to the reversal polynomial $\lambda^2K+\lambda D+M$ of $S(\lambda)$.

To see (iii) observe that by Theorem~\ref{thm:singind} the left minimal indices of $P(\lambda)$
are all zero and the right minimal indices of $P(\lambda)$ are at most one. By \cite[Theorem 5.10]{DeDM09a}
the left minimal indices of $S(\lambda)$ coincide with those of $P(\lambda)$, and if $\ve_1,\dots,\ve_k$
are the right minimal indices of $S(\lambda)$, then $\ve_1+1,\dots,\ve_k+1$ are the right minimal indices
of $P(\lambda)$. This implies that all minimal indices of $S(\lambda)$ are zero (if there are any).
\eproof
As another application of Theorem \ref{thm:singind} we obtain a result for rectangular matrix pencils which
in the case $E=I$, $Q>0$ coincides with the well known Lyapunov stability condition~\cite{LanT85} and can thus be seen as a generalization of this stability criterion.
In the following $Q^\dagger$ denotes the Moore-Penrose generalized inverse of a matrix $Q$, see, e.g.,~\cite{GolV96}.
\begin{corollary}\label{cor:ly1}
Let $A,E\in\Field^{n,m}$ and assume that there exists $Q\in\Field^{n,m}$ such that all minimal
indices of the pencil $\lambda E-Q$ are zero (if there are any) and
\begin{equation}
\label{eq:conditions}
E^\star Q\geq 0,\quad A Q^\dagger +Q^{\dagger\star} A^\star\leq 0,\quad \ker Q\subseteq\ker A.
\end{equation}
Then, the statements \mbox{\rm (i)--(v)} of Theorem~\mbox{\rm\ref{thm:singind}} hold for the pencil
$P(\lambda)=\lambda E-A\in\Field^{n,m}[\lambda]$.
\end{corollary}
\proof
We set $L:=AQ^\dagger$. The first and second inequality of \eqref{eq:conditions} ensure that the
assumptions of Theorem \ref{thm:singind} are satisfied for this $E,L,Q$. Furthermore, as
$\ker Q\subseteq\ker A$ and $Q^\dagger Q$ is the orthogonal projection onto $\im Q^\star\supset\im A^\star$,
we have $LQ=AQ^\dagger Q=A$.  Hence, the statements (i)--(v) of Theorem \ref{thm:singind} hold for the pencil $P(\lambda)$.
\eproof
In the case that $n=m$ and $Q$ is invertible, the condition $A Q^\dagger +Q^{\dagger\star} A^\star\leq 0$ can be
reformulated as $Q^\star A+A^\star Q\leq0$, and thus we obtain the following simplified statement.
\begin{corollary}\label{cor:ly2}
Let $A,E\in\Field^{n,n}$ and assume that there exist an invertible $Q\in\Field^{n,n}$
such that $E^\star Q\geq 0$ and $Q^\star A+A^\star Q\leq0$.
Then, the statements \mbox{\rm (i)--(v)} of Theorem~\mbox{\rm\ref{thm:singind}} hold for the pencil
$P(\lambda)=\lambda E-A\in\Field^{n,n}[\lambda]$ with the assumption of regularity of $\lambda E-Q$ in \mbox{\rm (v)} being automatically satisfied.
\end{corollary}

Other generalizations of Lyapunov's theorem with focus on asymptotic stability have been obtained in
\cite{IshT02,Sty02,TakMK95} for the case of square pencils and in \cite{IshT01} for the
rectangular case. Corollaries~\ref{cor:ly1} and~\ref{cor:ly2} give a generalization to
the concept of \emph{stability} in the sense of the conditions (i)--(v) in Theorem~\ref{thm:singind}.

\section{Proof of Theorem~\ref{thm:singind}}\label{sec:phproof}
In this section we will present a proof of our main result. To show the statements (iii)--(v)  we will need the following lemma.
\begin{lemma}\label{lem:1.11.17}
Let $E,A\in\mathbb F^{n,m}$ and $k\geq 2$. Assume that one of the following conditions holds
\begin{itemize}
\item[\rm (I)] the pencil $\lambda E-A$ has index $k$;
\item[\rm (II)] the pencil $\lambda E-A$ has a right minimal index $k-1$.
\end{itemize}
Then there exists a vector $x_1\in\mathbb F^{m}$ orthogonal to the common nullspace of $E$ and $A$ such that
the following two statements hold:
\begin{itemize}
\item[\rm (a)] There exist vectors $x_2,\dots,x_k\in\mathbb F^m$ such that $Ex_1=0$, $Ex_2=Ax_1$,
\dots, $Ex_k=Ax_{k-1}$, where in addition we have $Ax_{k}\neq 0$ in case {\rm (I)} and $Ax_k=0$ in case {\rm (II)}.
\item[\rm (b)] For any choice of vectors $x_2,\dots,x_k\in\mathbb F^m$ with
$Ex_1=0$, $Ex_2=A x_1$, \dots, $Ex_k=Ax_{k-1}$ we have
$Ax_1,\dots,Ax_{k-1}\neq 0$.
\end{itemize}
\end{lemma}
\proof
First, assume that $\lambda E-A$ is in Kronecker canonical form, where without loss of generality
the block $\mathcal N_k$ associated with a block of size $k$ of the eigenvalue $\infty$, or the block
associated with the right minimal index $k-1$, respectively, comes first, i.e., we have
\[
\lambda E-A=\mat{cc}\lambda E_{11}-A_{11}&0\\ 0&\lambda E_{22}-A_{22}\rix,
\]
where $\lambda E_{11}-A_{11}=\mathcal N_k$, or
\begin{equation}\label{eq:2.11.17}
\lambda E_{11}-A_{11}=\left[\begin{array}{cccc}
0&1\\&\ddots&\ddots\\&&0&1
\end{array}\right]-\lambda\left[\begin{array}{cccc}
1&0\\&\ddots&\ddots\\&&1&0
\end{array}\right].
\end{equation}
Here, in~\eqref{eq:2.11.17} we did on purpose not use the block $\mathcal L_{k-1}$ as the representation of a
block with right minimal index $k-1$, but its reversal which is well-known to be equivalent to $\mathcal L_{k-1}$.
In this way, we can treat the proofs for both cases simultaneously, because the block
in~\eqref{eq:2.11.17} consists exactly of the first $k-1$ rows of the block $\mathcal N_k$.

In the following, let $e_i^{(p)}$ denote the $i$-th standard basis vector of $\mathbb F^p$.
Then $x_1:=e_1^{(m)}$ is orthogonal to the common nullspace of $E$ and $A$, and
it is straightforward to check that the vectors $x_2=e_2^{(m)}$, \dots, $x_k=e_k^{(m)}$ satisfy (a).
On the other hand, if $x_2,\dots,x_k\in\mathbb F^m$ are chosen such that they satisfy
$Ex_1=0$, $Ex_2=Ax_1$, \dots, $E x_k=A x_{k-1}$ then a straightforward
computation shows that $x_2,\dots,x_k$ must have the form
\[
x_2=\mat{c}e_2^{(k)}+\alpha_1e_1^{(k)}\\ x_{22}\rix\ ,\dots,\
x_k=\mat{c}e_k^{(k)}+\alpha_1e_{k-1}^{(k)}+\dots+\alpha_{k-1}e_1^{(k)}\\ x_{k,2}\rix
\]
for appropriate $\alpha_1,\dots,\alpha_{k-1}\in\mathbb F$ and $x_{22},\dots,x_{k,2}\in\mathbb F^{m-k}$. Then we obtain
\[
Ex_j=Ax_{j-1}=\mat{c}e_{j-1}^{(k-1)}+\alpha_1e_{j-2}^{(k-1)}+\dots+\alpha_{j-2}e_1^{(k-1)}\\ A_{22}x_{j-1,2}\rix\neq 0
\]
for $j=2,\dots,k$ which shows (b).
This also proves (a) and (b) for the general case with the possible exception of the orthogonality condition on $x_1$. However, note that any vector $y\in\mathbb F^m$ satisfying $Ey=0=Ay$ can be
added to $x_1$, $x_2$, \dots $x_k$ without changing any of the identities
in (a) and (b). This shows that $x_1$ satisfying (a) and (b) can be chosen to be orthogonal
to the common nullspace of $E$ and $A$.
\eproof

\medskip

\noindent
\emph{Proof of Theorem~\ref{thm:singind}}:
Since $Q^\star E\geq 0$, we may assume without loss of generality by Theorem~\ref{thm:condensed}
that $E$ and $Q$ have the form
\begin{equation}\label{eq:30.9.17}
E=\mat{cc}E_{11}&0\\ 0&0\rix,\quad Q=\mat{cc}Q_{11}&0\\ 0&0\rix,
\end{equation}
respectively,  where $E_{11},Q_{11}$ are real and diagonal (which implies, in particular, that $Q^\star =Q)$
and satisfy $E_{11}^2+Q_{11}^2=I_{n_1}$. The fact that $E_{22}$ and $Q_{22}$ (and
thus also $E_{12}$ and $Q_{12}$) in~\eqref{cform} are zero follows from the assumption
that $\lambda E-Q$ has only left minimal indices equal to zero (if any).
To show the items (i)--(v) for the pencil $P(\lambda)=\lambda E-LQ$ we will again frequently use the decomposition $L=J-R$, where $R:=\frac{1}{2}(L+L^\star)\geq 0$  and $J:=\frac{1}{2}(L-L^\star)$.

(i) Let $\lambda_0\in\mathbb C$ be an eigenvalue of $P(\lambda)$ and let $v\neq 0$ be a regular eigenvector associated with
$\lambda_0$. Then we have $\lambda_0 Ev=LQv=(J-R)Qv$ and thus
\[
\lambda_0 v^\star Q^\star  Ev=v^\star Q^\star JQv-v^\star Q ^\star RQv.
\]
Considering the real parts of both sides of this equation, we obtain
\[
\operatorname{Re}(\lambda_0)\cdot v^\star Q^\star Ev=-v^\star Q ^\star RQv,
\]
where we used the fact that $Q^\star E$ and $R$ are Hermitian and $J$ is skew-Hermitian.
If $Q^\star Ev=0$, then the special structure of $E$ and $Q$ implies that $v=y_1+y_2$, where $Ey_1=0$ and
$Qy_2=0$ and thus $v$ can be expressed as a linear combination of vectors from the kernels of $E$ and $LQ$, respectively which by Lemma~\ref{lem:regdefsub} leads to a contradiction.
Hence, we have $Q^\star  Ev\neq 0$, and since $Q^\star  E$ is positive semidefinite, we obtain
$v^* Q^\star  Ev>0$ which finally implies
\[
\operatorname{Re}(\lambda_0)=-\frac{v^* Q^\star  RQv}{v^* Q^\star  Ev}\leq 0.
\]

For (ii)
we first prove the 'moreover' part.
For this, let the columns of $V\in\mathbb C^{m,k}$ form a basis of a regular
deflating subspace of $P(\lambda)$ associated with the eigenvalue $\lambda_0=i\omega$. Then we have to show $RQV=0$.

By the definition of regular deflating subspaces and from the Kronecker canonical form it follows
that there exists a matrix $W\in\mathbb C^{m,k}$ with full column rank such that
\begin{equation}\label{eq:10.10.17}
EV=W\quad\mbox{and}\quad (J-R)QV=WT,
\end{equation}
where $T\in\mathbb C^{k,k}$ only has the eigenvalue $i\omega$. Without
loss of generality we may assume that $T=i\omega I_k+N$ is in Jordan canonical form (JCF),
where $N$ is strictly upper triangular. Otherwise we choose $P$ such that $P^{-1}TP$ is in JCF
and set $\widetilde W=WP$ and $\widetilde V=VP$.

By the second identity in~\eqref{eq:10.10.17} we have $V^\star Q^\star (J-R)QV=V^\star Q^\star WT$,
and taking the Hermitian part of both sides we obtain
\begin{equation}\label{eq:27.10.17}
0\geq -2V^\star Q^\star  RQV=V^\star Q^\star WT+T^\star W^\star Q^\star V.
\end{equation}
Since $R$ is positive semidefinite, it remains to show that $V^\star Q^\star  RQV=0$, because then we also have $RQV=0$. For this, we first note that we have
\begin{equation}\label{VQEV}
V^\star Q^\star W=W^\star Q^\star V> 0.
\end{equation}
Indeed, it follows from the first identity in~\eqref{eq:10.10.17} that $V^\star Q^\star W= V^\star Q^\star  EV\geq 0$.
If there exists $x\neq 0$ such that $Q^\star EVx=0$, then with $y=Vx$ one has $Q Ey=0$. Due to the specific form of
$Q$ and $E$ this implies that $y=y_1+y_2$ with $Ey_1=0$ and $Qy_2=0$. Hence, $y$ is a vector belonging to a linear span of vectors from the kernels of $E$ and $LQ$ which, by Lemma~\ref{lem:regdefsub}, contradicts the fact that the columns of $V$
span a regular deflating subspace associated with $\lambda_0\neq 0$ and finishes the proof of \eqref{VQEV}.

Now let $M$ be the inverse of the Hermitian positive definite square root of the Hermitian positive definite matrix $V^\star Q^\star W=W^\star Q^\star V$. Then it follows from~\eqref{eq:27.10.17} that the matrix
\[
M(V^\star Q^\star WT+T^\star V^\star QW)M=M^{-1}TM+MT^\star M^{-1}
\]
is Hermitian negative semidefinite, because it is congruent to the right hand side of~\eqref{eq:27.10.17}.
Moreover, we obtain
\begin{eqnarray*}
\operatorname{trace}(M^{-1}TM+MT^\star M^{-1})&=&
\operatorname{trace}(M^{-1}TM)+\operatorname{trace}(MT^* M^{-1})\\
&=&\operatorname{trace}(T+T^\star)=\operatorname{trace}(N+N^\star)=0,
\end{eqnarray*}
because $N$ has a zero diagonal. But this implies that $M^{-1}TM+MT^\star M^{-1}=0$ and hence by~\eqref{eq:27.10.17}
also $-2V^\star Q RQ V=0$, which finishes the proof of the 'moreover' part.

Next, we will show that $i\omega\neq 0$ is a semisimple eigenvalue. For this, it remains to show that the matrix
$T=i\omega I_k+N$ in~\eqref{eq:10.10.17} is diagonal, i.e., $N=0$. Observe that with $RQV=0$ the
identities in~\eqref{eq:10.10.17} simplify to $EV=W$ and $JQV=WT$ which implies that
\[
V^\star Q^\star JQV=V^\star Q^\star WT.
\]
Using again the inverse $M$ of the Hermitian positive definite square root of $V^\star Q^\star WT$, we obtain that
\[
M^{-1}TM=MV^\star Q^\star WTM=MV^\star Q^\star JQVM,
\]
which implies that $T$ is similar to a matrix which is congruent to $J$, i.e., $T$ is similar to
a skew-Hermitian matrix. This, however, immediately implies that we must have $N=0$. Thus, $i\omega$  is a semisimple eigenvalue of  $\lambda E-(J-R)Q$ and assertion (ii) is proved.

For the remaining part of the proof, we will use a slight variation of~\eqref{eq:30.9.17}. By scaling the
nonzero parts of $E$ and $Q$ in~\eqref{eq:30.9.17} if necessary, we may assume without loss of generality that $E$ and $Q$ take the forms
\begin{equation}\label{eq:30.9.17a}
E=\mat{ccccc}I_k&0&0&0\\ 0&I_\ell&0&0\\ 0&0&0&0\\ 0&0&0&0\rix\quad\mbox{and}\quad
Q=\mat{ccccc}Q_1&0&0&0\\ 0&0&0&0\\ 0&0&I_q&0\\ 0&0&0&0\rix,
\end{equation}
where $Q_1\in\mathbb F^{k,k}$ is positive definite and both $E$ and $Q$ are partitioned conformably with
the partition $(k,\ell,q,p)$ of $n$, where $p,q\geq 0$.

To prove (iii), 
assume that the index $\nu$ of $P(\lambda)$ exceeds two. Then by Lemma~\ref{lem:1.11.17}
there exists a vector $x_1$ orthogonal to the common kernel of $E$ and $LQ$ that satisfies the
conditions (a) and (b) of Lemma~\ref{lem:1.11.17}. In particular, by (a) there exists vectors
$x_2,\dots,x_\nu$ such that
\begin{equation}\label{eq:30.9.17b}
Ex_1=0,\; Ex_2=LQx_1\neq 0,\; Ex_3=LQx_2\neq 0,\;\dots,\; Ex_\nu=LQx_{\nu-1}\neq 0,\; LQx_\nu\neq 0.
\end{equation}
Since $x_1$ is in the kernel of $E$, it must be of the form
$x_1^\top=\mat{cccc}0&0&x_{13}^\top&0\rix^\top$, when  partitioned conformably with respect to the partition $(k,\ell,q,p)$ of $n$. Then, we can apply a simultaneous
similarity transformation to $E$, $L$, and $Q$ with a unitary transformation matrix of the form
$U=I_k\oplus I_\ell\oplus U_q\oplus I_{p}$, where $U_mx_{13}=\alpha e_q$ for some $\alpha\in\mathbb R$, e.g., we can take $U_q$ to be an appropriate Householder matrix. Note that the simultaneous similarity
transformation with $U$ preserves all the properties of $E$, $Q$, and $L$. Finally, $x_1$ can be appropriately scaled such that $\alpha=1$.

Next, let us extend the partitioning of $E$, $Q$, $L$, $J$, and $R$ by splitting the third block row and column
in two, i.e., the new partitioning is according to the partition $(k,\ell,q-1,1,p)$
of $n$. Then we have
\[
E=\operatorname{diag}(I_k,I_\ell,0,0,0),\quad Q=\operatorname{diag}(Q_1,0,I_{q-1},1,0)
\]
and we can express $L=J-R$ as
\[
L=\mat{ccccc}L_{11}&L_{12}&L_{13}&J_{14}-R_{14}&L_{15}\\
L_{21}&L_{22}&L_{23}&J_{24}-R_{24}&L_{25}\\
L_{31}&L_{32}&L_{33}&J_{34}-R_{34}&L_{35}\\
-J_{14}^\star-R_{14}^\star&L_{42}&L_{43}&J_{44}-R_{44}&L_{45}\\
L_{51}&L_{52}&L_{53}&-J_{45}^\star-R_{45}^\star&L_{55}\rix.
\]
 Observe that
\[
LQx_1=\mat{c}J_{14}-R_{14}\\ J_{24}-R_{24}\\ J_{34}-R_{34}\\ J_{44}-R_{44}\\-J_{45}^\star-R_{45}^\star\rix\in\Field^{n}
\]
has to be in the range of $E$ in order to satisfy $Ex_2=LQx_1$ for some $x_2$. This implies the identities
\begin{equation}\label{eq:30.9.17c}
J_{34}=R_{34},\quad J_{44}=R_{44},\quad -J_{45}^\star=R_{45}^\star.
\end{equation}
{F}rom this we obtain $J_{44}=R_{44}=0$, because $R_{44}$ is a real scalar and $J_{44}$ is purely imaginary.
But then it follows from the positive-semidefiniteness of $R$ that $R_{14}=0$, $R_{24}=0$, $R_{34}=0$, and
$R_{45}=0$ which in turn implies that $J_{34}=0$ and $J_{45}=0$ by~\eqref{eq:30.9.17c}. Thus, using
$L_{43}=-J_{34}^\star+R_{34}^\star=0$ we find that $LQ$ and
$LQx_1$ take the simplified form
\[
LQ=\mat{ccccc}L_{11}Q_1&0&L_{13}&J_{14}&0\\
L_{21}Q_1&0&L_{23}&J_{24}&0\\
L_{31}Q_1&0&L_{33}&0&0\\
-J_{14}^\star Q_1&0&0&0&0\\
L_{51}Q_1&0&L_{53}&0&0\rix\quad\mbox{and}\quad
LQx_1=\mat{c}J_{14}\\ J_{24}\\ 0\\ 0\\0\rix,
\]
where we used the fact that the second diagonal block in $Q$ is zero.
At this point, we see that $x_2$ has to be of the form $x_2^\top=\mat{ccccc}J_{14}^\top&J_{24}^\top&
x_{23}^\top&x_{24}^\top&x_{25}^\top\rix^\top$ in order to satisfy
\[
\mat{c}J_{14}\\ J_{24}\\ 0\\ 0\\0\rix=LQx_1=Ex_2=\mat{ccccc}I&0&0&0&0\\ 0&I&0&0&0\\ 0&0&0&0&0\\ 0&0&0&0&0\\ 0&0&0&0&0\rix
x_2.
\]
Since  $Ex_3=LQx_2$, the vector
\[
LQx_2=\mat{c}\ast\\\ast\\\ast\\ -J_{14}^\star Q_1 J_{14}\\ \ast\rix
\]
has to be in the range of $E$, which implies that $J_{14}^\star Q_1J_{14}=0$, and hence $J_{14}=0$ as $Q_1$ is positive definite. But then the vector
$\widetilde x_2:=\mat{ccccc}0&J_{24}^\top&0&0&0\rix^\top$ satisfies the identities
\[
Ex_1=0,\quad E\widetilde x_2=LQx_1,\quad LQ\widetilde x_2=0,
\]
which is a contradiction to part (b) of Lemma~\ref{lem:1.11.17}

To prove (iv)
assume that $P(\lambda)$ has a right minimal index $\varepsilon\geq 2$. Then, by Lemma~\ref{lem:1.11.17},
there exists a vector $x_1$, orthogonal to the common nullspace of $E$ and $LQ$, that satisfies the conditions (a) and (b) of Lemma~\ref{lem:1.11.17}. In particular, by (a) there exists vectors satisfying
\begin{equation*}\label{eq:30.9.17}
Ex_1=0,\; Ex_2=LQx_1\neq 0,\; Ex_3=LQx_2\neq 0,\;\dots,\; Ex_{\varepsilon+1}=LQx_{\varepsilon}\neq 0,
\; LQx_{\varepsilon+1}=0.
\end{equation*}
Since only $x_1,x_2,x_3$ have been used in
(iii) to obtain a contradiction, we can follow exactly the lines of (iii) to obtain again a contradiction.

(v)
Since $\lambda E-Q$ is regular, we have $p=0$ in the partitioning~\eqref{eq:30.9.17a}.
Assume that $P(\lambda)$ has a left minimal index of size $\eta>0$. Then by applying
Lemma~\ref{lem:1.11.17} to the pencil $P(\lambda)^\star=\lambda E^{\star}-(LQ)^\star$, there exists a vector $x_1$, orthogonal to the common nullspace of $E$ and $LQ$, such that the conditions (a) and (b) of Lemma~\ref{lem:1.11.17} are satisfied. In particular, (a) implies the existence of vectors $x_2,\dots,x_{\eta+1}$ such that
\[
x_1^\star E=0,\quad x_2^\star E=x_1^\star LQ\neq 0,\quad\cdots,\quad x_{\eta+1}^\star E=x_\eta^\star LQ\neq 0,
\quad x_{\eta+1}^\star LQ=0.
\]
By an argument analogous to the one in the proof of (iii), we may assume
without loss of generality that $x_1=e_n$ is the last standard basis vector of $\Field^n$.
Next, partitioning $E$, $Q$, $J$, and $R$ conformably
with respect to the partition $(k,\ell,q-1,1)$ of $n$ we get
\begin{eqnarray*}\nonumber 
E&=&\mat{cccc}I_k&0&0&0\\ 0&I_\ell&0&0\\ 0&0&0&0\\ 0&0&0&0\rix,\
Q=\mat{cccc}Q_1&0&0&0\\ 0&0&0&0\\ 0&0&I_{q-1}&0\\ 0&0&0&1\rix,\\\label{29.9.17two}
J&=&\mat{cccc}J_{11}&J_{12}&J_{13}&J_{14}\\ -J_{12}^\star & J_{22}& J_{23}&J_{24}\\
-J_{13}^\star & -J_{23}^\star & J_{33}&J_{34}\\ -J_{14}^*&-J_{24}^*&-J_{34}^*&J_{44}\rix,
R=\mat{cccc}R_{11}&R_{12}&R_{13}&R_{14}\\ R_{12}^\star & R_{22}& R_{23}&R_{24}\\
R_{13}^\star & R_{23}^\star & R_{33}&R_{34}\\ R_{14}^*&R_{24}^*&R_{34}^*&R_{44}\rix.
\end{eqnarray*}
Then, similarly as in the proof of (iii), we obtain the identity
\[
x_2^*E=x_1^*LQ=\mat{cccc} (-J_{14}^\star-R_{14}^\star)Q_1&0&-J_{34}^\star-R_{34}^\star
&J_{44}-R_{44}\rix
\]
which implies that $J_{44}=R_{44}=0$, $J_{34}=R_{34}=0$, $R_{24}=0$, and $R_{14}=0$. Hence,  the matrix $LQ$ takes the simplified form
\begin{equation*}\label{29.9.17lq}
LQ=(J-R)Q=
\mat{cccc}(J_{11}-R_{11})Q_1&0&J_{13}-R_{13}&J_{14}\\ (-J_{12}^\star-R_{12}^\star)Q_1 & 0& J_{23}-R_{23}&J_{24}\\
(-J_{13}^\star-R_{13}^\star)Q_1 & 0 & J_{33}-R_{33}&0\\ -J_{14}^\star Q_1&0&0&0\rix.
\end{equation*}
This implies that the vector $x_2$ must be of the form
$x_2^\star=\mat{cccc}-J_{14}^\star Q_1& 0&x_{23}^\star& x_{24}^\star\rix$
for some $x_{23}\in\mathbb F^{q-1},x_{24}\in\mathbb F$ in order to satisfy the identity
\begin{equation}\label{eq:xLQsim}
x_2^\star E=x_1^\star LQ=\mat{cccc} -J_{14}^\star Q_1 & 0 & 0 & 0\rix.
\end{equation}
Let us assume first that $\eta>1$. Then the vector
\[
x_2^\star LQ=\mat{cccc}\ast&\ast&\ast&-J_{14}^\star Q_1J_{14}\rix
\]
has to satisfy the identity $x^\star_3 E=x_2^\star LQ$, which implies
$J_{14}^\star Q_1J_{14}=0$ and thus $J_{14}=0$ as $Q_1$ is positive definite.
But this together with \eqref{eq:xLQsim} implies that $x_1^\star LQ=0$ which contradicts the fact that $x_1$ is orthogonal to the common left nullspace of $E$ and $LQ$.

In the second case that $\eta=1$ we have a chain $x_1^\star E=0$, $x_2^\star E=x_1^\star LQ\neq 0$
and $x_2^\star LQ=0$, which  again leads to $J_{14}=0$ and $x_1^\star LQ= 0$, a contradiction.
Thus, $P(\lambda)$ cannot have a left minimal index of size larger than zero which finishes the proof.
\eproof

\section{Removing higher order Jordan blocks of the eigenvalue $0$}\label{sec5}
In this section we will concentrate on the case that the pencil
\begin{equation}\label{phpencilind1}
P(\lambda)=\lambda E-LQ, \quad E,Q,L\in\mathbb F^{n,n}
\end{equation}
with $Q^\star E=E^\star Q\geq 0$ and $L+L^\star\leq 0$ is regular.
As we have seen in the previous section, it may happen that the index of the pencil is two.
In \cite{BeaMXZ17} it was shown that in the case of port-Hamiltonian DAEs with variable coefficients
it is possible to perform an index reduction of the given port-Hamiltonian DAE such that the resulting DAE of index one is again a port-Hamiltonian DAE. Since the case with constant coefficients is a special
case of the port-Hamiltonian DAEs considered in \cite{BeaMXZ17}, the index reduction performed in Section~6 of that paper can also be applied in our case, so that we can assume that our pencil $P(\lambda)$ is of index at most one. However, as it was shown in the previous section, it is still possible that the eigenvalue $0$ of $P(\lambda)$ is not semisimple. In this section, we will therefore provide a method for checking this property
this and we also present a method of perturbing $L$ to make $0$ a semisimple eigenvalue, while keeping the structure assumptions on the pencil and all other eigenvalues and
their Kronecker structures intact. The main ingredient of the method is a condensed form that provides particular information on the Kronecker structure of the eigenvalue zero.
\begin{theorem}\label{zerosemi}
Assume that $P(\lambda)$ from \eqref{phpencilind1} with $E^\star Q=Q^\star E\geq 0$ and $L+L^\star\leq 0$
is regular and of index at most one.
Then there exist a partition $(n_1,n_2,n_3,n_4)$ of $n$, a unitary matrix $\widetilde U$ and an invertible matrix
$\widetilde X$ such that partitioning conformably with respect
to $(n_1,n_2,n_3,n_4)$ we have
\begin{eqnarray} \nonumber 
\widetilde Q&:=&\widetilde UQ \widetilde X=\matp{ Q_{11} &  Q_{12} &0 & 0\\
 Q_{21} &  Q_{22} &0 & 0\\  0 & 0 &0 & 0\\ Q_{41} & Q_{42} & 0 & I},\
\widetilde E:=\widetilde UE \widetilde X=\matp{ E_{11} &0 &0 & 0\\
 E_{21} &  E_{22} &0 & 0\\  0 & 0 &I & 0\\ 0 & 0 & 0 & 0},\\
\label{Bdeco2}
\widetilde L&:=&\widetilde U L\widetilde U^\star=\mat{c} L_{ij}\rix_{ij=1}^4,\ \widetilde A:=\widetilde ULQ \widetilde X=\matp{A_{11} &0 &0 & A_{14}\\ A_{21} & 0 &0 & A_{24}\\
A_{31} & A_{32} &0 & A_{34}\\ 0 & 0 & 0 & A_{44}},
\end{eqnarray}
 where  $E_{11}$, $ E_{22} $, $A_{11}$ are lower-triangular and invertible, and $A_{44}$,
$\matp{Q_{11} & Q_{12} \\ Q_{21} & Q_{22}}$ are invertible.
\end{theorem}
\proof
Note that the regularity of $P(\lambda)$ implies the regularity of $\lambda E -Q$ by Proposition~\ref{prop:regular}.
Hence, by Proposition \ref{prop:EQreg}, there exist a unitary matrix $U$ and invertible matrix $X$ such that
\[
\hat Q:=UQX=\matp{Q_1 &&\\ & 0 & \\ && I},\quad
\hat E:=UEX=\matp{E_1 &&\\ & I & \\ && 0 },
\]
with $E_1$ and $Q_1$ being  diagonal and positive. The sizes of the block columns and rows
are $n_1+n_2,n_3,n_4$ in both matrices, later on the first block row and column will be split into blocks of sizes $n_1$, $n_2$, respectively.
Next, let
\[
\hat L:=ULU^*=\matp{\hat L_{ij}}_{ij=1,\dots,3},\quad ULQX=
 \matp{\hat L_{11} Q_1 &0  & \hat L_{13} \\ \hat L_{21}Q_1 &0 & \hat L_{22} \\ \hat L_{31}Q_1 & 0 &\hat L_{33}}.
\]
be partitioned conformably with $\hat Q$ and $\hat E$. Since $P(\lambda)$ is of index one,
it follows that $\hat L_{33}$ is invertible, see, e.g., \cite{BreCP96}. Setting
\begin{equation}\label{eq:T}
 T=\matp{I &0&0\\ 0& I &0\\ -\hat L_{33}^{-1}\hat L_{31}Q_1 & 0& I}
\end{equation}
we get that
\begin{equation}\label{eq:ULQXT}
UEXT=\hat E,\quad U\hat LQXT=\matp{(\hat L_{11}-\hat L_{13}\hat L_{33}^{-1}\hat L_{31}) Q_1 &0  & \hat L_{13} \\
(\hat L_{21}-\hat L_{23}\hat L_{33}^{-1}\hat L_{31})Q_1 &0 & \hat L_{23} \\0  & 0 &\hat L_{33}}.
\end{equation}
Observe that
\begin{equation}\label{Lupleft}
Q_1^\star E_1^{-1}> 0,\quad (\hat L_{11}-\hat L_{13}\hat L_{33}^{-1}\hat L_{31})+(\hat L_{11}-\hat L_{13}\hat L_{33}^{-1}\hat L_{31})^\star \leq 0 .
\end{equation}
Indeed, the first inequality results directly from the form of $Q_1$ and $E_1$ and the second follows from
$\hat L+\hat L^\star\leq 0$ as
\[
\hat L_{11}-\hat L_{13}\hat L_{33}^{-1}\hat L_{31}=S^\star \hat  L \hat S  ,\quad S:=\mat{c} I \\ 0 \\ -\hat L_{33}^{-1}\hat L_{31}\rix.
\]
Hence, we can apply \cite[Lemma 3.1]{MehMS16} to see that zero is a semisimple eigenvalue of the pencil
$
\lambda I -(\hat L_{11}-\hat L_{13}\hat L_{33}^{-1}\hat L_{31})Q_1 E_1^{-1}
$
and thus also a semisimple eigenvalue of the pencil
$\lambda E_1- (\hat L_{11}-\hat L_{13}\hat L_{33}^{-1}\hat L_{31})Q_1 $, which has a generalized Schur form,
see, e.g.,~\cite{GolV96},
\begin{equation}\label{eq:W}
 W_1 (\hat L_{11}-\hat L_{13}\hat L_{33}^{-1}\hat L_{31}) Q_1 W_2=\matp{A_{11} & 0\\
             A_{21} & 0},\quad W_1E_1 W_2=\matp{ E_{11} & 0\\
              E_{21} &  E_{22}}
\end{equation}
with $A_{11},E_{11}\in\Field^{n_1,n_1}$, $E_{22}\in\Field^{n_2,n_2}$ being lower triangular and invertible and $W_1,W_2$ being unitary.
Setting
\[
\widetilde U=\matp{W_1 \\&I\\ && I}U,\quad \widetilde X= XT\matp{W_2\\ &I \\&& I},
\]
and splitting the first column and row conformably with $(n_1,n_2)$, where
$n_1=\dim\ker(\hat L_{21}-\hat L_{23}\hat L_{33}^{-1}\hat L_{31})Q_1$, finishes the definition of
$\tilde Q,\tilde X$ and $(n_1,n_2,n_3,n_4)$.

The particular form of $\tilde A$ and $\tilde E$ follows now from \eqref{eq:T}, \eqref{eq:ULQXT} and \eqref{eq:W}. Furthermore,
\[\tilde Q= \matp{W_1 \\&I\\ && I}\hat Q T\matp{W_2\\ &I \\&& I}
\]
and hence
$
\matp{Q_{11} & Q_{12} \\ Q_{21} & Q_{22}}=W_1Q_1W_2
$
is invertible and the form of $\tilde Q$ follows.
\eproof
Let us discuss now some properties of the condensed form \eqref{Bdeco2}.
\begin{corollary}\label{cor:semi-eq}
With the same notation and assumptions as in Theorem~\mbox{\rm\ref{zerosemi}},
the following conditions are equivalent:
\begin{itemize}
\item[\rm (a)] Zero is a semisimple eigenvalue of $P(\lambda)$,
\item[\rm (b)] $A_{32}=0$,
\item[\rm (c)]  $\im  \mat{c} Q_{21} \\ Q_{22}\\ 0\\ Q_{24}\rix\subseteq \ker \mat{cccc} L_{31} & L_{32} &0 & L_{34}\rix$.
\end{itemize}
\end{corollary}
\proof (a)$\Leftrightarrow$(b):
As $A_{44}$ is invertible, the pencil $\lambda\widetilde E-\widetilde A$ is equivalent to the pencil
\begin{equation*}\label{Bdeco1}
\lambda\matp{E_{11} &0 &0 & 0\\ E_{12} & E_{22} &0 & 0\\  0 & 0 &I & 0\\ 0 & 0 & 0 & 0}-
\matp{A_{11} &0 &0 & 0\\ A_{21} & 0 &0 & 0\\  A_{31} & A_{32} &0 & 0\\ 0& 0 & 0 & A_{44}}.
\end{equation*}
which has a semisimple eigenvalue at $0$ if and only if the matrix
\[
C:=\matp{E_{11} &0 &0 \\ E_{12} & E_{22} &0\\  0 & 0 &I}^{-1} \matp{A_{11} &0 &0 \\ A_{21} & 0 &0 \\
 A_{31} & A_{32} &0}=\matp{E^{-1}_{11}A_{11} & 0 &0\\ * & 0 &0 \\  A_{31} & A_{32} &0 },
\]
has a semisimple eigenvalue at 0. Since $E_{11}^{-1}A_{11}$ is invertible,
we immediately see that semi-simplicity of the eigenvalue zero is equivalent to $A_{32}=0$.

(b)$\Leftrightarrow$(c):
Let $L_j$ (respectively, $Q_j$)  denote for $j=1,2,3,4$ the $j$-th block row (block column) of the matrix $\widetilde L$ (matrix $\widetilde Q$).
Then $A_{32}=L_3Q_2$ and (c) is clearly equivalent to $L_3Q_2=0$.
%
%
\eproof
Following the notation of \eqref{Bdeco2} we also obtain the following result.
\begin{corollary}\label{cor:M}
With the same notation and assumptions as in Theorem~\mbox{\rm\ref{zerosemi}},
any perturbation $L+M$ of $L$, where  $\widetilde M:=\widetilde U M \widetilde U^\star$ is of the form
\[
\widetilde M=\matp{0 & 0 & M_{13} & 0\\ 0 & 0 & M_{23} &0 \\ M_{31} & M_{32} & M_{33} & M_{34} \\ 0 & 0 & M_{43} & 0}  ,
\]
leaves the pencil $P(\lambda)$ regular and the Kronecker form invariant, except for possible changes in
the number and sizes of the blocks corresponding to the eigenvalue $0$.
\end{corollary}
\proof
The assertion follows directly by computing
\[
\widetilde A+ \widetilde M\widetilde Q=\matp{A_{11} &0 &0 & A_{14}\\ A_{21} & 0 &0 & A_{24}\\
A_{31}+W_{31} & A_{32}+W_{32} &0 & A_{34}+W_{34}\\ 0 & 0 & 0 & A_{44}},\quad W=\mat{c} W_{ij}\rix _{i,j=1}^4:=\widetilde M\widetilde Q.
\]
Due to the form of $\widetilde E$ in \eqref{Bdeco2}, the entries $W_{31}$ and $W_{34}$ do not have any
influence on the Kronecker form and regularity of $\lambda \widetilde E - \widetilde A$ (and thus also
$\lambda E-A$) and the entry $W_{32}$ may cause changes to the Kronecker structure of the eigenvalue zero only.
\eproof
Hence, given a pencil $P(\lambda)$ as in~\eqref{phpencilind1} with an eigenvalue $0$ that is not semisimple, it is desirable to construct a stabilizing perturbation of the form above with $W_{32}=-A_{32}$ as this procedure will make the eigenvalue zero semisimple, but leaves all other eigenvalues and their Kronecker structure invariant. To express such a perturbation in terms of $L$ and $Q$ we need some preparations. In the following, for a linear subspace $\mathcal X\subseteq\Field ^n$ the symbol
$P_{\mathcal X}\in\Field^{n,n}$  denotes the orthogonal projection onto $ \mathcal X$ and
$\mathcal X\ominus \mathcal Y$ stands for $P_{{\mathcal Y}^\perp}\mathcal X$.
\begin{lemma}\label{ABZ}
Let $C\in\Field^{m,l}$, $B\in\Field^{n,m}$ and let $Z_0\in\Field^{m,l}$ be such that $B(C+Z_0)=0$. Then
$\im C\ominus\ker B\subseteq \im Z_0\subseteq \im C\oplus\ker B$ and the matrix $Z$ with the smallest Frobenius
norm satisfying $B(C+Z)=0$ has the property that $\im C\ominus\ker B=\im Z$ and its Frobenius norm
equals $\norm{Z}_F=\norm{P_{\im B^*}C}_F$.
\end{lemma}
\proof
The proof follows directly from the decomposition $\Field^m=\ker B\oplus\im B^\star$.
\eproof
Theorem 1.1 of \cite{FoiF90} gives a characterization of positive semidefinite $2\times2$ block
operators. We present here a matrix version for the reader's convenience.
\begin{theorem}\label{l1}
Let $B\in\mathbb F^{k,k}$ and $F\in\Field^{m,m}$ be Hermitian positive semidefinite, let $E\in\mathbb F^{m,k}$ and let
\[
T:=\matp{B & E^\star \\ E & F   }.
\]
Then $T\geq 0$ if and only if there exist a linear mapping $\Gamma:\im F\to\im B$ with $\norm{\Gamma}\leq 1$ such that $E^\star=B^{1/2}\Gamma F^{1/2}$. In particular, if $T\geq 0$ then $\im E^\star\subseteq\im B^{1/2}=\im B$.
\end{theorem}
%
%
The following  lemma will serve as a prerequisite for constructing perturbations of $R$ that
preserve  semi-definiteness.
\begin{lemma}\label{l2}
Let $B\in\mathbb F^{k,k}$ and $D\in\Field^{m,m}$ be Hermitian, positive semidefinite, and let $C\in\mathbb F^{m,k}$ be such that
\[
\matp{B & C^\star \\ C & D   }\geq 0.
\]
 Then, for any $Y\in\Field^{m,k}$ with $\im Y^\star \subseteq \im B$
we have that
\[
T(Y,W):=\matp{B & (C+Y)^\star \\ C+Y & D+ W }\geq 0,
\]
for $W:=(\alpha^2+2\alpha)( D+\norm Y I)+\norm Y I\;$ with
$\;\alpha=\big \|(B\vert_{\im B})^{-\frac12} Y^\star (D+\norm Y I)^{-\frac12}\big \|$.
\end{lemma}
\proof
  Theorem \ref{l1} adapted to the current situation states that $T(Y,W)$ is positive semidefinite if and only if there exists a linear map $\Gamma:\im(D+W)\to\im B$ with $\norm{\Gamma}\leq 1$ such that
\begin{equation}\label{Gammaeq}
(C+Y)^\star=B^{1/2}\ \Gamma\  (D+W)^{1/2}.
\end{equation}
Let now $W$ be  as in the assertion and assume that $Y\neq 0$. The case $Y=0$ is trivial, because then $W=0$ as well. Observe that $D\geq 0$ implies that $D+W=(\alpha+1)^2(D+\norm YI)>0$. Setting
\[
\Gamma:=(B \vert_{\im B})^{-\frac12}(  C+Y )^\star (D+W)^{-\frac12}=(1+\alpha)^{-1}(B \vert_{\im B})^{-\frac12}( C+Y )^\star (D+\norm Y I)^{-\frac12},
\]
by Lemma~\ref{l1} we have that $\im C^*\subseteq\im B$ and by the assumption that $\im Y^*\subseteq\im B$ we have that \eqref{Gammaeq} is satisfied. Furthermore, by definition of $\alpha$, we have that
\[
\norm{\Gamma}\leq(1+\alpha)^{-1} \left(\big \|(B\vert_{\im B})^{-\frac12} C^\star (D+\norm Y I)^{-\frac12}\big \| +\alpha\right).
\]
Note that by Theorem \ref{l1}, applied to $T(0,\norm Y I)\geq 0$, there exist $\tilde\Gamma:\Field^m\to\im B$
with $\|\tilde \Gamma\|\leq 1$ such that $C^\star=B^{1/2}\ \tilde \Gamma\  (D+\norm Y I)^{1/2}$.
Hence,
\[
\big \|(B\vert_{\im B})^{-\frac12} C^\star (D+\norm Y I)^{-\frac12}\big \|\leq\|\tilde \Gamma\|\leq 1,
\]
%
which finishes the proof.
\eproof
We have now prepared the tools for  constructing a stabilizing perturbation of a
pencil as in~\eqref{phpencilind1} by perturbing the matrix $L$ in a structure preserving manner. Below $\norm{X}_F$ stands for the Frobenius norm of the matrix $X$.
\begin{theorem}\label{th:delta}
Assume that $P(\lambda)$ from \eqref{phpencilind1}, with $E^\star Q=Q^\star E\geq 0$ and $L+L^\star\leq 0$,
is a regular pencil of index at most one, and let the decomposition \eqref{Bdeco2} be given. Then, for every $Y:=\matp{Y_1 & Y_2 & Y_4}$ with
\[
\im Y^\star\subseteq \im R_0, \quad R_0:=\matp{ R_{11} & R_{12}  & R_{14}\\  R_{21} & R_{22} &  R_{24}\\
 R_{41} & R_{42}  & R_{44}  } ,
\]
there exist
 $\Delta_R=\Delta_R^\star$ such that $R+\Delta R\geq 0$ and
\begin{equation}\label{eq:DeltaR}
 \norm{\Delta_R}\leq 2\norm Y +\left(\big\| (R_0\vert_{\im R_0})^{-1} \big\|\norm Y+ 2\big\| (R_0\vert_{\im R_0})^{-\frac12} \big\| \norm{ Y}^{\frac12}\right)(\norm{R_{33}}+\norm Y)
\end{equation}
and there exist $\Delta_J=-\Delta_J^\star$ with the Frobenius norm  bounded by
\begin{equation}\label{eq:DeltaJ}
\norm{\Delta_J}_F\leq 2\norm{P_{\im Q_2}( L_3^\star+Y^\star )  }_F,
\end{equation}
where
\[
L_3=\matp{L_{31} & L_{32} & L_{33} & L_{34}},\ Q_2=\matp{Q_{21}\\ Q_{22}\\ 0 \\ Q_{24}},
\]
 such that for every $s,t\in\Comp$ the pencil
$\lambda E+ (J+s\Delta J- R - t \Delta_R)Q$ is regular and has the same Kronecker structure for the nonzero eigenvalues (including infinity)
as $P(\lambda)$ and zero is a semisimple eigenvalue of $\lambda E- (J+\Delta_J-R-\Delta_R)Q$.
\end{theorem}
\proof
Without loss of generality we may assume that $\widetilde U=I_n$ in the decomposition
\eqref{Bdeco2}. We interchange the last two block rows and last two block columns, so that
\[
R=\matp{R_0  & C^* \\ C & R_{33}  },\quad R_3:=\matp{R_{31} & R_{32}  & R_{34}}.
\]
Let $Y=\matp{Y_1 & Y_2 & Y_4}$ with $\im Y^\star \subseteq\im R_0$ and let $B:=R_0$, $D=R_{33}$.
In this setting we apply Lemma \ref{l2} getting a matrix $W$ such that $R+\Delta R\geq 0$, where
\[
\Delta_R:=\matp{0 & Y^\star \\ Y& W} .
\]
The bound \eqref{eq:DeltaR} then follows, since
\[
\big\| (R_{33}+\norm Y I)^{-\frac12} \big\| \leq \big\| (\norm Y I)^{-\frac12}\big\| =\norm{Y}^{-\frac 12}
\]
and hence $
\alpha \leq \big\| (R_0\vert_{\im R_0})^{-\frac12} \big\| \norm{ Y}^{\frac12}$,
with $\alpha$ as in Lemma~\ref{l2}.
Furthermore, by Corollary \ref{cor:M} for any $t\in\Comp$ the pencil $\lambda E- (J- R - t \Delta_R)Q$
is regular and has the same Kronecker structure associated with the nonzero eigenvalues (including infinity).

Now let us define $\Delta_J$ to make the eigenvalue zero of the perturbed pencil semisimple. We reverse the interchange of the last two block rows and block columns, so that we are in the original setting of \eqref{Bdeco2}, and we define $Z=\matp{Z_1 & Z_2 &0& Z_4}$ as in  Lemma~\ref{ABZ} applied
to $B=Q_2^\star$, $C=L_3^\star +\hat Y^\star$, where
$L_3$ and $Q_2$ are defined as in the assertion and $\hat Y=\matp{Y_1 & Y_2 & 0 & Y_4}$.
Then we get $Q_2^\star( L_3^\star +\hat Y^\star + Z^\star  )=0  $, and defining
\[
\Delta_J:=\matp{0 & 0 & -Z_1^\star & 0\\     0&0 &-Z_2^\star &0\\Z_1 & Z_2 &0& Z_4\\ 0 & 0 & -Z_4^\star & 0},
\]
in the perturbed matrix $A+\Delta_R+\Delta_J$ the $(3,2)$ block entry equals zero, which
by Corollary \ref{cor:semi-eq} is equivalent to semi-simplicity of the zero eigenvalue. The estimate  \eqref{eq:DeltaJ} follows by Lemma~\ref{ABZ} and again by Corollary~\ref{cor:M}
for any $s\in\Comp$ the pencil $\lambda E-(J+s\Delta_J - R - t \Delta_R)Q$ is regular and has
the same Kronecker structure at nonzero eigenvalues (including infinity).
\eproof

Let us conclude  this section with  some comments.
%
\begin{itemize}
\item Theorem~\ref{th:delta} can be applied with $Y=0$. In this way we  get a skew-symmetric perturbation of $L$.
\item Note that condition (c) of Corollary \ref{cor:semi-eq} is implied by
\[
  \ker \mat{cccc} L_{11} & L_{12}& 0 & L_{14}\\  L_{21} & L_{22}& 0 & L_{24}\rix\subseteq \ker \mat{cccc} L_{31} & L_{32}& 0 &L_{34}\rix,
\]
and this leads to a bound
$\norm{\Delta_J}_F\leq 2\norm{P_{\im [L_1,L_2]}( L_3^\star+Y^\star )  }_F$,
which, in general, is weaker than \eqref{eq:DeltaJ}, however, it does not use the matrix $Q$.
\item By Theorem~\ref{l1}, if $\im Y^\star$ is not contained in $R_0$ then the matrix $R+\Delta_R$ is not positive semidefinite.
\item For sufficiently small $Y$ with $\im ( Y+Y^\star)\subseteq \im R$ the perturbation
$\widetilde\Delta R:=Y+Y^\star$ preserves positive definiteness as well. Note that
$\|\widetilde\Delta_R\|=\norm{Y}$, while $\norm{\Delta_R}=\mathcal{O}(\norm{Y}^{1/2})$ as
$\norm Y$ tends to zero. However, note that the condition $\im ( Y+Y^\star)\subseteq \im R$ is stronger than $\im Y^\star\subseteq \im R_0$ as it involves $R_{33}$ as well.
Furthermore, in applications, $Y$ does not need to have a small norm.
\item The matrix $\Delta_R$ constructed in the proof of Theorem~\ref{th:delta} has the same  form as the matrix $M$ from Corollary~\ref{cor:M}. Hence, $\Delta_R$ is never positive semidefinite (unless it is zero). Therefore, strengthening the statement from  $R+\Delta_R$ positive semidefinite to $\Delta_R$
positive semi-definite is in general not possible.
\item The case $\Delta_J=0$ can be considered as well. A sufficient condition for the existence of a symmetric perturbation $\Delta_L$ (of the form presented in the theorem) is that $P_{\im Q_2}L_3^\star\subseteq \im R_0$. In this case  the matrix $Y=-L_3P_{\im Q_2}$ satisfies $\im Y^\star\subseteq\im R_0$ and for this particular choice of $Y$ we get $\Delta J=0$.
\item The matrices $\widetilde U$ and $\widetilde X$, as well as the decomposition \eqref{Bdeco2} can be computed efficiently, e.g., using the CS-decomposition and the generalized Schur decomposition of \texttt{matlab}.
\end{itemize}

\section*{Conclusions}
We have studied matrix pencils associated with dissipative Hamiltonian
descriptor systems. We have derived canonical and condensed forms and characterized spectral properties of such pencils. In particular, we have studied  the size of Jordan blocks associated with the eigenvalues $0$ and $\infty$ and the restrictions for the left and right minimal indices.
When Kronecker blocks of size larger than one arise for the eigenvalue $0$ (i.e. the associated system is not Lyapunov stable) then we have constructed minimal norm structure-preserving perturbations that make the system Lyapunov stable.
\bibliographystyle{plain}

\end{document}